\documentclass[11pt]{article}

\usepackage[T1]{fontenc}
\usepackage{lmodern}
\usepackage[margin=1in]{geometry}
\usepackage{amsmath,amssymb,amsthm,mathtools}
\usepackage{microtype}
\usepackage{booktabs,tabularx}
\usepackage{pdfpages}
\usepackage{tikz}
\usetikzlibrary{arrows.meta,calc,decorations.pathreplacing,positioning}
\usepackage[hidelinks]{hyperref}

\newtheorem{theorem}{Theorem}[section]
\newtheorem{proposition}[theorem]{Proposition}
\newtheorem{lemma}[theorem]{Lemma}
\newtheorem{corollary}[theorem]{Corollary}
\newcommand{\restatedname}{}
\newtheorem*{restatedthm}{\restatedname}
\newenvironment{restated}[1]
 {\renewcommand{\restatedname}{Theorem \ref{#1}}\restatedthm}
 {\endrestatedthm}
\newenvironment{restatedL}[1]
 {\renewcommand{\restatedname}{Lemma \ref{#1}}\restatedthm}
 {\endrestatedthm}
\theoremstyle{definition}
\newtheorem{definition}[theorem]{Definition}
\newtheorem{example}[theorem]{Example}

\newtheorem*{yabloexample}{Example}

\newcommand{\post}[1]{N_{#1}}
\newcommand{\pre}[1]{N^-_{#1}}
\newcommand{\feeds}{\rightsquigarrow}

\newcommand{\bd}{\partial}
\newcommand{\two}{\{0,1\}}
\tikzset{>=Stealth,
 dot/.style={circle,fill,inner sep=1.3pt},
 odot/.style={circle,draw,fill=white,inner sep=1.3pt},
 gate/.style={rectangle,fill,inner sep=2pt},
 box/.style={draw,rounded corners=3pt,align=center,inner sep=4pt,font=\small}}

\newcommand{\noo}[1]{}

\title{Parity-safe digraphs with countably many ends\\ are kernel-perfect}
\author{Sjur Dyrkolbotn, Claude Fable 5.2, Micha{\l} Walicki}
\date{22.11.2026}
\hypersetup{pdftitle={Parity-safe digraphs with countably many ends are kernel-perfect}}

\begin{document}
%%%% ---- part-A
\maketitle

\begin{abstract}
A kernel of a digraph is an independent set that every other vertex
points into, and a digraph is kernel-perfect if every induced
subgraph has a kernel.
A digraph without odd cycles in which no
ray contains infinitely many vertices that each send an infinite fan
of disjoint paths to that ray is \emph{safe}. 
The conjecture from \cite{b10} that every safe
digraph is kernel-perfect was proved there for finitely many ends.

We prove a stronger statement for at most countably many ends.  Number the vertices of a
ray, and let the \emph{type} of a path that ends on the ray be the parity of its length
plus the number of its last vertex. A fan is \emph{mixed} if it has infinitely many paths
of each type, and a digraph is \emph{parity-safe} if it has no odd cycle and no ray
containing infinitely many vertices with a mixed fan to it. Yablo's graph $(\omega,<)$ is
not parity-safe. Its variant with chords only between levels of different parity has a
kernel -- it is parity-safe but not safe. Our main result is that \emph{a parity-safe
  digraph with at most countably many ends is kernel-perfect}, and the claim for safe
digraphs follows since safety implies parity-safety.  The digraph may have any
cardinality, and nothing is assumed about the structure of the ends, but the case of
uncountably many ends remains open.
\end{abstract}

\section{Introduction}
\label{s1}

\subsection{Kernels and paradox}
\label{s1.1}

A digraph $\Gamma$ has a vertex set $V(\Gamma)$ and edges $u\to v$, for $u,v\in V(\Gamma)$. A set $K$
of vertices is a \emph{kernel} if it contains no edge and every
vertex outside $K$ has an edge into $K$. Think of each vertex as a
sentence saying ``everything I point to is false'':
\[
 p_v\ \leftrightarrow\ \bigwedge_{v\to w}\neg p_w .
\]
A truth assignment satisfies all of these equations exactly when the
true vertices form a kernel. A digraph without a kernel is a
paradox.\footnote{This name comes from the application of kernels to paradox analysis, e.g., \cite{Cook,JPL}.} A directed cycle of odd length is one, but a cycle of even length is not; it
has two kernels.

\begin{yabloexample}[Yablo]
\emph{Vertices $v_0,v_1,\dots$ with $v_i\to v_j$ for all $i<j$. If $v_i$
were true, all later vertices would be false; but then $v_{i+1}$,
pointing only to false vertices, would be true. So every vertex is
false, and then $v_0$ should be true. There is no kernel, although
there is no cycle.}
\begin{center}
\begin{tikzpicture}[x=1.5cm]
 \foreach \i in {0,...,4} \node[dot,label=below:{\small$v_{\i}$}] (v\i) at (\i,0){};
 \node at (5,0){$\cdots$};
 \foreach \i/\j in {0/1,1/2,2/3,3/4} \draw[->] (v\i)--(v\j);
 \foreach \i/\j in {0/2,1/3,2/4} \draw[->] (v\i) to[bend left=35] (v\j);
 \foreach \i/\j in {0/3,1/4} \draw[->] (v\i) to[bend left=50] (v\j);
 \draw[->] (v0) to[bend left=62] (v4);
\end{tikzpicture}
\end{center}
\end{yabloexample}

A digraph is \emph{kernel-perfect} if every induced subgraph has a
kernel. This is the robust form of ``contains no paradox''. By
Richardson's theorem \cite{b9} a finitely branching digraph without odd cycles
is kernel-perfect. So without odd cycles paradox needs infinity, and
the question is which infinite configurations are responsible.

\subsection{Fans and parity fans}
\label{s1.2}

A \emph{ray} is an infinite path $r_0\to r_1\to\cdots$ with distinct
vertices. A vertex $v$ \emph{dominates} a ray $R$ if it has an infinite fan to $R$, namely, 
infinitely many paths from $v$ to $R$ that pairwise share only $v$.

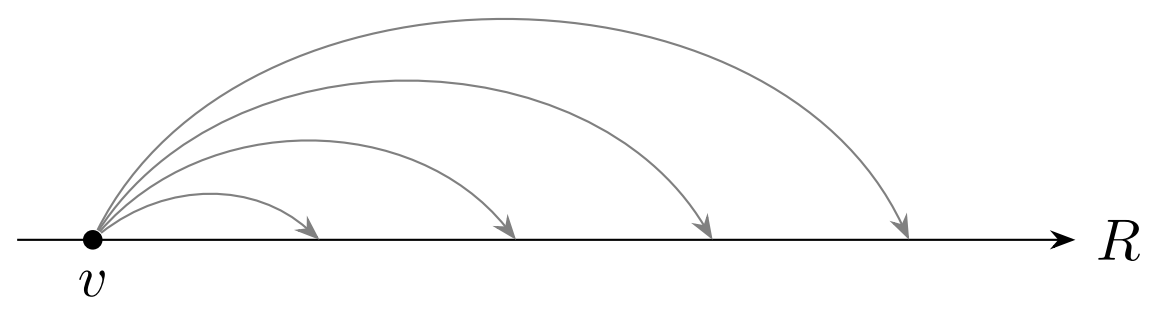
\begin{center}
\begin{tikzpicture}
 \draw[->] (0,0)--(7,0) node[right]{$R$};
 \node[dot,label=below:{$v$}] (v) at (.5,0){};
 \foreach \x/\b in {2/40,3.3/50,4.6/58,5.9/64}
   \draw[->,gray] (v) to[bend left=\b] (\x,0);
\end{tikzpicture}
\end{center}

In Yablo's graph every vertex dominates the ray $v_0\to v_1\to\cdots$ on which it lies. In
\cite{b10} a digraph is called \emph{safe} if it has no odd cycle and no ray containing
infinitely many vertices that dominate it. It is conjectured that \emph{every safe digraph
  is kernel-perfect}, \cite[Conjecture 1.2]{b10}).

This is there proven for digraphs with finitely many ends
\cite[Theorem 3.19]{b10} and for some classes with infinitely
many ends of a special form. An \emph{end} is, roughly, a direction in
which rays go. 
A ray $R$ \emph{feeds} a ray $R'$ if $R$ keeps reaching $R'$
past every finite obstacle. Mutual feeding is an equivalence relation and an \emph{end} is such an equivalence class of rays.

Paradox is a matter of parity and therefore odd cycles must be excluded. From now on, a
\emph{digraph} is arbitrary and is denoted by $\Gamma$, while a \emph{noc graph} is a
digraph with no odd cycles and is denoted by $G$.  Loops are odd cycles, so noc graphs
have none, but may have even cycles and any cardinality.

Safety considers fans but pays no attention to parity.  In Yablo's graph the trouble is
that $v_i$ reaches $v_j$ both in one step and in two.

\begin{example}[Odd Yablo]
\label{n1.1}
Keep only the edges $v_i\to v_j$ of Yablo's graph with $j-i$ odd
(the removed edges are dotted).
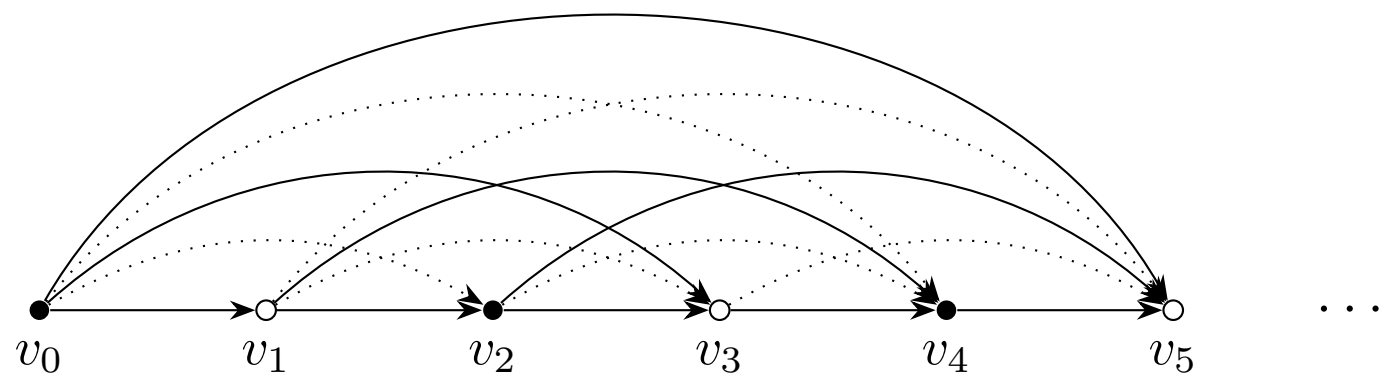
\begin{center}
\begin{tikzpicture}[x=1.5cm]
 \foreach \i in {0,2,4} \node[dot,label=below:{\small$v_{\i}$}] (v\i) at (\i,0){};
 \foreach \i in {1,3,5} \node[odot,label=below:{\small$v_{\i}$}] (v\i) at (\i,0){};
 \node at (5.8,0){$\cdots$};
 % removed edges: $j-i$ even
 \foreach \i/\j in {0/2,1/3,2/4,3/5} \draw[->,dotted] (v\i) to[bend left=30] (v\j);
 \foreach \i/\j in {0/4,1/5} \draw[->,dotted] (v\i) to[bend left=52] (v\j);
 % remaining edges: $j-i$ odd
 \foreach \i/\j in {0/1,1/2,2/3,3/4,4/5} \draw[->] (v\i)--(v\j);
 \foreach \i/\j in {0/3,1/4,2/5} \draw[->] (v\i) to[bend left=42] (v\j);
 \draw[->] (v0) to[bend left=60] (v5);
\end{tikzpicture}
\end{center}
Every vertex still dominates the ray $v_0\to v_1\to\cdots$ by
direct edges, so this digraph is not safe. But every edge joins an even
index to an odd one. The digraph is bipartite, the vertices with even
index form a kernel, and it is kernel-perfect.
\end{example}
Thus, safety forbids too much in terms of fans. What Yablo's graph has, and the odd
version lacks, is a fan whose paths disagree about parity. Let
$R=r_0\to r_1\to\cdots$ be a ray and call $n$ the \emph{level} of
$r_n$. The \emph{type} of a path $P$ that ends at $r_n$ is the
parity of $|P|+n$, where $|P|$ is the length of $P$. An infinite fan
from $v$ to $R$ is \emph{mixed} if it has infinitely many paths of
each type. If there is one, $v$ \emph{parity-dominates} $R$. A
digraph is \emph{parity-safe} if it has no odd cycle and no ray
containing infinitely many vertices that parity-dominate it. Safe
digraphs are parity-safe.

In Yablo's graph the edges from $v_i$ to the later vertices form a
mixed fan. In the odd version all paths from $v_i$ to the ray have
the same type, so no fan is mixed. For a vertex $v=r_m$ on the ray, a
path from $v$ to $r_n$ has the type of $m$ exactly if it agrees in
parity with the segment of $R$ from $v$ to $r_n$. So a mixed fan
from a vertex of the ray contains infinitely many disjoint paths that
\emph{disagree} in parity with the ray: odd connections between
levels of equal parity, or even ones between levels of different
parity.

\subsection{The result and the structure of the proof}\label{s1.3}

\begin{theorem}\label{n1.2}
A parity-safe digraph with at most countably many ends is
kernel-perfect. In particular, a safe digraph with at most countably
many ends is kernel-perfect.
\end{theorem}

The digraph may have any cardinality, and nothing is assumed about the structure of its
ends.  The heart of the proof is Theorem \ref{n6.4} by which a parity-safe digraph that is
not kernel-perfect has uncountably many \emph{ladder ends}, defined as follows:

\begin{itemize}
\item a \emph{finite-port region} is a set $X$ of vertices with finitely many out-neighbours (outside $X$); 
\item a vertex is \emph{essential} if no finite-port region containing 
it induces a kernel-perfect digraph; 
\item a \emph{ladder ray} is a ray through infinitely many essential
vertices $b_1,b_2,\dots$ such that from each $b_i$ to each later
$b_j$ there are paths of both parities;
\item a \emph{ladder end} is an end containing a ladder ray.
\end{itemize}
The proof of Theorem \ref{n6.4} rests on the following principle, which
may be of independent interest. 

\begin{lemma}[Finite-port cover]
\label{n1.3}
If every vertex of a noc graph $G$ lies
in a finite-port region that induces a kernel-perfect digraph, then
$G$ is kernel-perfect.
\end{lemma}
The regions may overlap in any way, and there may be any number of
them. The lemma fails for digraphs with odd cycles: in a directed
triangle every single vertex is a kernel-perfect finite-port region.
If every vertex has finitely many successors, each single vertex is
a kernel-perfect finite-port region; so the lemma contains
Richardson's theorem. % \cite{b9}.

Overall, the proof of Theorem \ref{n1.2} is a diagonalisation argument showing that if a parity-safe digraph $\Gamma$, with countably many ends
$e_0,e_1,e_2,\dots$, is not kernel-perfect, then it has an end not on the list. The main steps are as follows:

\begin{enumerate}
\item 
A region with finitely many ports
interacts with the rest of the digraph through finitely many truth
values. If every vertex had a paradox-free region of this kind
around it, the regions could be glued and $\Gamma$ would be
kernel-perfect (Lemma \ref{n1.3}). So some vertices are
essential, and inside every paradoxical region there is a ray
through infinitely many essential vertices
(Proposition \ref{n5.5}).
\item 
Let $L$ be such a ray. By
parity-safety some essential vertex $q$ on $L$ sends no mixed fan to
$L$. Then there is a finite set $F$, the gate, such that all paths
from $q$ to a tail of $L$ that avoid $F$ have the same type.
\item 
Let $Z^*$ be everything $q$ can reach without passing the gate. It is
a finite-port region, so, $q$ being essential, it contains a
paradox. Split $Z^*$ into the part $A$ that can still reach the tail
of $L$ and the rest, $Z'$. Because all paths from $q$ through $A$ to
$L$ have the same type, $A$ is bipartite, like the odd Yablo graph,
and bipartite digraphs are harmless. So the paradox sits in $Z'$,
which therefore contains another ray $L'$ of the same kind. No ray
inside $Z'$ can keep reaching $L$.
\item 
Do this for $e_0$, then inside the new
region for $e_1$, and so on. Walk along $L_0$ to $q_0$, cross into
$Z_1$, walk along $L_1$ to $q_1$, cross into $Z_2$, and so on. The
resulting ray $R$ has a tail in every $Z_i$, so its end is none of
the $e_i$ (Theorem \ref{n6.4}).
\end{enumerate}

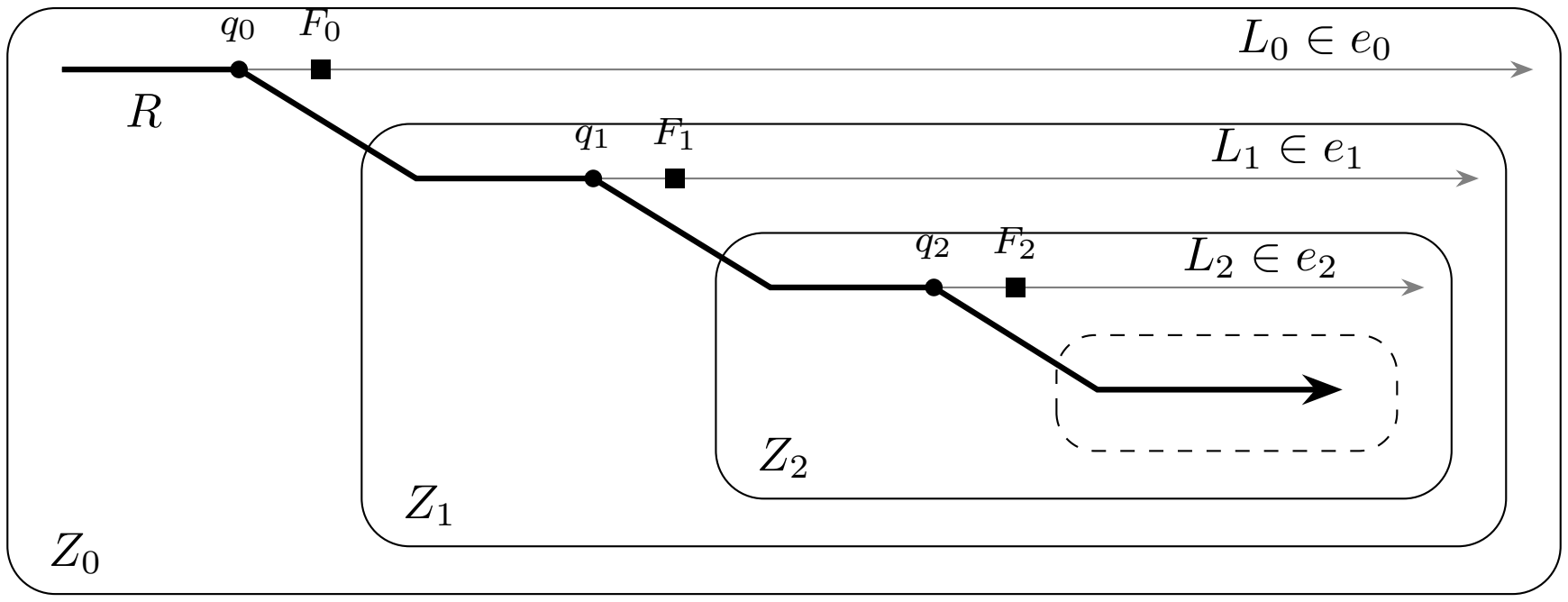
\begin{center}
\begin{tikzpicture}
 \draw[rounded corners=10pt] (0,0) rectangle (11.4,4.3);
 \draw[rounded corners=10pt] (2.6,.35) rectangle (11,3.45);
 \draw[rounded corners=10pt] (5.2,.7) rectangle (10.6,2.65);
 \draw[rounded corners=8pt,dashed] (7.7,1.05) rectangle (10.2,1.9);
 \node at (.5,.3){\small$Z_0$};
 \node at (3.1,.65){\small$Z_1$};
 \node at (5.7,1){\small$Z_2$};
 \draw[->,gray] (.4,3.85)--(11.2,3.85);
 \node at (9.6,4.07){\small$L_0\in e_0$};
 \draw[->,gray] (3,3.05)--(10.8,3.05);
 \node at (9.4,3.27){\small$L_1\in e_1$};
 \draw[->,gray] (5.6,2.25)--(10.4,2.25);
 \node at (9.2,2.47){\small$L_2\in e_2$};
 \node[gate,label=above:{\scriptsize$F_0$}] at (2.3,3.85){};
 \node[gate,label=above:{\scriptsize$F_1$}] at (4.9,3.05){};
 \node[gate,label=above:{\scriptsize$F_2$}] at (7.4,2.25){};
 \node[dot,label=above:{\scriptsize$q_0$}] at (1.7,3.85){};
 \node[dot,label=above:{\scriptsize$q_1$}] at (4.3,3.05){};
 \node[dot,label=above:{\scriptsize$q_2$}] at (6.8,2.25){};
 \draw[->,very thick] (.4,3.85)--(1.7,3.85)--(3,3.05)--(4.3,3.05)--(5.6,2.25)
   --(6.8,2.25)--(8,1.5)--(9.8,1.5);
 \node at (1,3.55){\small$R$};
\end{tikzpicture}

\smallskip
{\small The region $Z_{i+1}$ lies behind the gate $F_i$, cannot reach
the tail of $L_i$, and contains no ray of the end $e_i$. The thick
ray $R$ has a tail in every $Z_i$.}
\end{center}

\subsection*{The structure of the paper and the use of Claude}
\label{s1.6}

Section \ref{s2} presents the standard notions and results and Section \ref{s3} some basic
concepts specific for this paper. Section \ref{s4} introduces finite-port regions and
proves the central Lemma \ref{n1.3}. Section \ref{s5} defines essential vertices and
ladder rays. Section \ref{s6} presents the diagonal argument and shows that a
non-kernel-perfect parity-safe digraph has uncountably many ladder ends, from which
Theorem \ref{n1.2} follows. A reader interested only in acyclic digraphs can ignore every
sentence about strong components: they are then single vertices.

The results were developed through several hours of interaction with Claude, during which
the tool was guided and carried out most of the formalisation. The results have been
verified in Lean, and Section \ref{s8} presents the Lean formalisation used for the main
Theorem \ref{n1.2}, together with the necessary definitions. The proof itself can be found at \url{https://walickimichal.github.io/michal.github.io/count-parity-safe-proof.lean}.

\section{Kernels and semikernels}\label{s2}
$\Gamma[X]$ is the subgraph induced by $X$, and $\Gamma-F=\Gamma[V(\Gamma)\setminus
F]$. Write $\post{\Gamma}(u)$ for the set of successors (out-neighbours) of $u$ in
$\Gamma$ and $\pre{\Gamma}(u)$ for its predecessors in $\Gamma$, and extend both to sets
of vertices. For a subgraph $H$, $\post{H}(u)$ is the set of successors of $u$ in $H$;
thus $\post{\Gamma[X]}(u)=\post{\Gamma}(u)\cap X$. A \emph{walk} is a finite sequence of
vertices in which each points to the next; its \emph{length} is the number of its edges,
and it is \emph{closed} if it ends where it starts. A \emph{path} is a walk with distinct
vertices; a single vertex is a path of length zero. A \emph{cycle} is a closed walk of
positive length that repeats no vertex except the first. A vertex $u$ \emph{reaches} $v$
if there is a path from $u$ to $v$; every vertex reaches itself. A set $D$ is
\emph{successor-closed} if $\post{\Gamma}(D)\subseteq D$. The \emph{future} of a vertex is
the set of vertices it reaches. A \emph{sink} is a vertex without successors. The
\emph{strong component} of $v$ is the set of vertices $w$ such that $v$ reaches $w$ and
$w$ reaches $v$.

A set is \emph{independent} if it contains no edge. A \emph{kernel}
is an independent set $K$ such that every vertex outside $K$ has a
successor in $K$. In terms of the truth assignment $x$ that is $1$
on $K$ and $0$ elsewhere, $K$ is a kernel if and only if the
\emph{kernel equation}
\[
 x(v)=1\quad\Longleftrightarrow\quad x(w)=0\ \text{for all }w\in\post{\Gamma}(v)
\]
holds at every vertex $v$.

\begin{lemma}[Parity principle]
\label{n2.1}
Let $G$ be a noc graph.
\begin{enumerate}
\item Every closed walk has even length.
\item Every walk from $u$ to $v$ contains a path from $u$ to $v$
whose length has the same parity. Here ``contains'' means that the
path uses only vertices of the walk.
\item Every strong component $S$ has a colouring $c:S\to\two$ such
that every edge inside $S$ joins vertices of different colours. All
walks between two given vertices of $S$ have the same parity.
\end{enumerate}
\end{lemma}
\begin{proof}
(1) Take a closed walk of odd length that is as short as possible.
If it repeated a vertex other than its first, it would split into
two shorter closed walks, one of them odd. So it is an odd cycle.

(2) If the walk repeats a vertex, cut out the closed walk between
two occurrences. By (1) the parity does not change. Repeat.

(3) Fix $r\in S$. If two walks from $r$ to $x\in S$ had different
parities, then together with a walk from $x$ to $r$ they would give
closed walks of both parities. So the parity $c(x)$ of the walks
from $r$ to $x$ is well defined, and it changes along every edge
inside $S$. A walk between two vertices of $S$ stays in $S$, so its
parity is determined by the colours of its ends.
\end{proof}

A digraph is \emph{bipartite} if all its vertices can be coloured
with two colours so that every edge joins different colours. It then
has no odd cycle. By (3) every strong component of a noc graph induces a bipartite
digraph. We fix such a colouring $c$ of every strong component once
and for all.

A \emph{semikernel} is an independent set $S$ with
$\post{\Gamma}(S)\subseteq\pre{\Gamma}(S)$: whatever $S$ points to points back into
$S$, \cite{NL71}. It is a partial solution that cannot be spoiled from outside.
The next three statements hold for all digraphs. They are essentially given in \cite{NL71} which  established also the subsequent proposition (using Zorn's lemma).

\begin{lemma}[\cite{NL71}]
\label{n2.2}
\begin{enumerate}
\item The union of a chain of semikernels is a semikernel.
\item If $D$ is successor-closed, every semikernel of $\Gamma[D]$, in
particular every kernel of $\Gamma[D]$, is a semikernel of $\Gamma$.
\item Let $S$ be a semikernel of $\Gamma$ and
$C=V(\Gamma)\setminus(S\cup\pre{\Gamma}(S))$. If $S'$ is a semikernel of $\Gamma[C]$,
then $S\cup S'$ is a semikernel of $\Gamma$.
\end{enumerate}
\end{lemma}
% \begin{proof}
% (1) An edge inside the union, or a successor of the union, already
% involves a single member of the chain. (2) Successors in $\Gamma$ of a
% subset of $D$ are its successors in $\Gamma[D]$. (3) There is no edge
% from $S$ to $C$, because $\post(S)\subseteq\pre(S)$, and none from
% $C$ to $S$ by the definition of $C$. So $S\cup S'$ is independent. A
% successor of $S$ points into $S$. A successor of $S'$ is not in $S$;
% if it lies in $C$ it points into $S'$, and otherwise it lies in
% $\pre(S)$ and points into $S$.
% \end{proof}

\begin{proposition}[\cite{NL71}]
\label{n2.3}
If every nonempty induced subgraph of $\Gamma$ has a nonempty semikernel,
then $\Gamma$ is kernel-perfect.
\end{proposition}
% \begin{proof}
% The hypothesis passes to induced subgraphs, so it suffices to find a
% kernel of $\Gamma$. By Lemma \ref{n2.2}(1) and Zorn's lemma there is a
% maximal semikernel $S$. Let $C=V(\Gamma)\setminus(S\cup\pre(S))$. If $C$
% were nonempty, $\Gamma[C]$ would have a nonempty semikernel $S'$, and
% $S\cup S'$ would be a larger semikernel. So $C=\varnothing$, and $S$
% is a kernel.
% \end{proof}
A \emph{core} is a nonempty digraph without a nonempty semikernel.

\begin{corollary}\label{n2.4}
A digraph that is not kernel-perfect contains an induced core. In a
core $C$ there is no sink, and for every nonempty successor-closed
$D\subseteq V(C)$ the digraph $C[D]$ is a core and has no kernel. In
particular the future of every vertex of $C$ induces a core.
\end{corollary}
\begin{proof}
The first statement is the contrapositive of
Proposition \ref{n2.3}. A sink is a one-element semikernel. For the
rest use Lemma \ref{n2.2}(2).
\end{proof}

\begin{lemma}[Two pieces]
\label{n2.5}
Let $D\subseteq V(\Gamma)$ be successor-closed. If $\Gamma[D]$ and $\Gamma-D$ are
kernel-perfect, then $\Gamma$ is kernel-perfect.
\end{lemma}
\begin{proof}
Let $H$ be a nonempty induced subgraph of $\Gamma$. If $H$ has no vertex
in $D$, it has a kernel. Otherwise $V(H)\cap D$ is successor-closed
in $H$ and induces a digraph with a kernel, which is a nonempty
semikernel of $H$ by Lemma \ref{n2.2}(2). Apply
Proposition \ref{n2.3}.
\end{proof}

\subsection{Canonical kernels of bipartite digraphs}\label{s2.4}
Bipartite digraphs are kernel-perfect. We need a little more: a
kernel that is chosen in a definite way and responds monotonically
when vertices are removed. Let $B$ be a bipartite digraph with colour
classes $B_0,B_1$, and let $Q\subseteq V(B)$. For $T\subseteq Q\cap
B_0$ put
\[
 \Phi_Q(T)=\{x\in Q\cap B_0:\ \text{every successor of $x$ in $Q$ has a successor in $T$}\}.
\]
Let $T(Q)$ be the union of all sets $T\subseteq Q\cap B_0$ with
$T\subseteq\Phi_Q(T)$, and
\[
 \kappa(Q)=T(Q)\cup\{y\in Q\cap B_1:\ \text{$y$ has no successor in $T(Q)$}\}.
\]
Read $T(Q)$ as the largest set of vertices of colour $0$ that can
consistently be declared true.

\begin{lemma}\label{n2.6}
\begin{enumerate}
\item $\kappa(Q)$ is a kernel of $B[Q]$. So every bipartite digraph
is kernel-perfect.
\item If $Q\cap B_0\subseteq Q'\cap B_0$ and $Q\cap B_1\supseteq
Q'\cap B_1$, then $\kappa(Q)\cap B_0\subseteq\kappa(Q')\cap B_0$ and
$\kappa(Q)\cap B_1\supseteq\kappa(Q')\cap B_1$.
\end{enumerate}
\end{lemma}
\begin{proof}
The map $\Phi_Q$ is monotone: a larger $T$ gives a larger
$\Phi_Q(T)$. Write $T=T(Q)$. Each set $T_0$ with
$T_0\subseteq\Phi_Q(T_0)$ satisfies $T_0\subseteq\Phi_Q(T)$, so
$T\subseteq\Phi_Q(T)$. Applying $\Phi_Q$ gives
$\Phi_Q(T)\subseteq\Phi_Q(\Phi_Q(T))$, so $\Phi_Q(T)$ is one of the
sets in the union, and $\Phi_Q(T)\subseteq T$. Hence
$T=\Phi_Q(T)$.

(1) Edges of $B[Q]$ join $B_0$ to $B_1$. Let $x\in T$ and
$y\in\kappa(Q)\cap B_1$. There is no edge $y\to x$ by the definition
of $\kappa$, and none $x\to y$, because $x\in\Phi_Q(T)$ means that
every successor of $x$ in $Q$ has a successor in $T$. So $\kappa(Q)$
is independent. A vertex $y\in Q\cap B_1$ outside $\kappa(Q)$ has a
successor in $T$. A vertex $x\in Q\cap B_0$ outside $T=\Phi_Q(T)$
has a successor $y\in Q$ without a successor Later in $T$; this $y$ lies
in $\kappa(Q)$.

(2) Let $T=T(Q)$. If $x\in T$, then $x\in Q'\cap B_0$, and every
successor of $x$ in $Q'$ lies in $Q$ and so has a successor in $T$.
Hence $T\subseteq\Phi_{Q'}(T)$ and $T\subseteq T(Q')$. A vertex of
$\kappa(Q')\cap B_1$ lies in $Q$ and has no successor in
$T(Q')\supseteq T$, so it lies in $\kappa(Q)$.
\end{proof}
For kernels of a single bipartite digraph, condition (2) can be reduced to only one part (intersection with only $B_0$ or only $B_1$), which gives a complete lattice, but later, in Lemma \ref{n4.8}, this order is applied in a way requiring both components. 
Exchanging the roles of $B_0$ and $B_1$ gives the other (greatest or least) kernel, which satisfies (2) as well; nothing below depends on this choice. 

\section{Rays: fans, parity fans, feeding and ends}
\label{s3}

A \emph{ray} is an infinite sequence $R=r_0\to r_1\to\cdots$ of
distinct vertices; its \emph{tails} are the rays
$r_n\to r_{n+1}\to\cdots$. For vertices of a ray, ``beyond'' and
``up to'' refer to the order along the ray. A ray \emph{lies in} a
set $Z$ if all its vertices are in $Z$.

An \emph{infinite fan} from a vertex $v$ to a ray $R$ is an infinite
family of paths of positive length from $v$ to vertices of $R$, any
two of which share only $v$. If there is one, $v$ \emph{dominates}
$R$. The vertex $v$ may or may not lie on $R$, and the paths may
meet $R$ before they end.

Let $R=r_0\to r_1\to\cdots$ be a ray. The \emph{level} of $r_n$ is
$n$. For a walk $P$ that ends at $r_n$, the \emph{type} of $P$ with
respect to $R$ is the parity of $|P|+n$, where $|P|$ is the length
of $P$.

\begin{definition}
\label{n3.2}
An infinite fan from $v$ to $R$ is \emph{mixed} if infinitely many
of its paths have even type and infinitely many have odd type. If
there is a mixed fan from $v$ to $R$, then $v$
\emph{parity-dominates} $R$. A digraph is \emph{parity-safe} if it has no odd cycle and no ray
contains infinitely many vertices that parity-dominate it.
\end{definition}

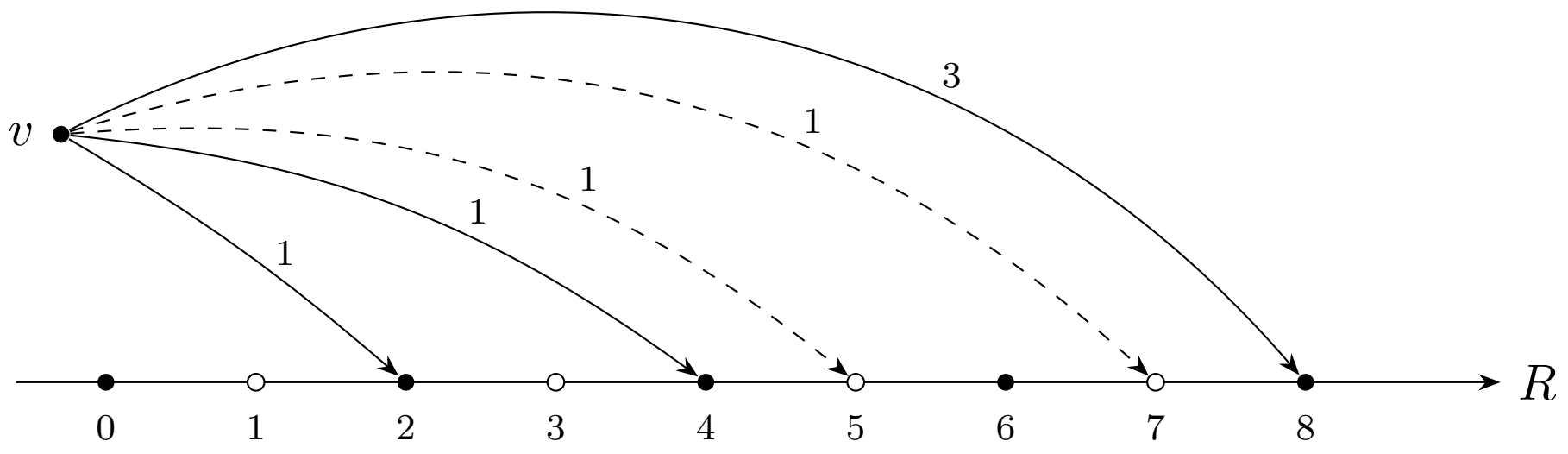
\begin{center}
\begin{tikzpicture}[x=1.15cm]
 \draw[->] (-.6,0)--(9.3,0) node[right]{$R$};
 \foreach \i in {0,2,4,6,8} \node[dot] (r\i) at (\i,0){};
 \foreach \i in {1,3,5,7} \node[odot] (r\i) at (\i,0){};
 \foreach \i in {0,...,8} \node[below=3.4pt] at (\i,0){\scriptsize$\i$};
 \node[dot,label=left:{$v$}] (v) at (-.3,1.9){};
 \draw[->] (v) to[bend left=5] node[pos=.65,above]{\scriptsize$1$} (r2);
 \draw[->] (v) to[bend left=15] node[pos=.63,above]{\scriptsize$1$} (r4);
 \draw[->,dashed] (v) to[bend left=22] node[pos=.64,above]{\scriptsize$1$} (r5);
 \draw[->,dashed] (v) to[bend left=30] node[pos=.66,above]{\scriptsize$1$} (r7);
 \draw[->] (v) to[bend left=38] node[pos=.68,above]{\scriptsize$3$} (r8);
\end{tikzpicture}

\smallskip
{\small Levels on the ray and lengths on the fan paths. Solid paths
have odd type (length plus level is odd), dashed paths even type. A
fan with infinitely many of each is mixed.}
\end{center}

Three remarks on the definition.
\begin{itemize}
\item Passing to a tail of $R$ lowers all levels by the same number.
This flips all types or none. So whether $v$ parity-dominates $R$
does not depend on where the numbering starts, and $v$
parity-dominates $R$ if and only if it parity-dominates a tail of
$R$.
\item Parity-domination implies domination. So \emph{safe digraphs are
parity-safe}. Parity-safety passes to induced subgraphs.
\item If $v=r_m$ lies on $R$, a walk $P$ from $v$ to $r_n$
\emph{agrees with the ray} if $|P|$ has the parity of $n-m$, and
\emph{disagrees} otherwise. Agreeing walks have the type of $m$,
disagreeing walks the other type.
\end{itemize}

\begin{lemma}[Parity cut criterion]
\label{n3.3}
A vertex $v$ parity-dominates a ray $R$ if and only if for every
finite set $F\not\ni v$ and every tail $T$ of $R$ there are two
paths of positive length in $\Gamma-F$ from $v$ to vertices of $T$ whose
types with respect to $R$ are different.

So if $v$ does \emph{not} parity-dominate $R$, there are a finite
set $F\not\ni v$ and a tail $T$ of $R$ such that all paths of
positive length in $\Gamma-F$ from $v$ to $T$ have the same type. The set
$F$ may be enlarged by any finite set not containing $v$, and $T$
may be replaced by any of its tails.
\end{lemma}
\begin{proof}
Let $v$ send a mixed fan to $R$. A finite set $F\not\ni v$ meets
only finitely many of its paths. Infinitely many paths of each type
remain, with distinct ends, so one of each type ends in $T$.
Conversely, we build a mixed fan $P_1,P_2,\dots$ in which $P_k$ has
the type of $k$. If $P_1,\dots,P_{k-1}$ have been found, let $F$
consist of their vertices other than $v$. There are two paths of
positive length in $\Gamma-F$ from $v$ to $R$ with different types. Let
$P_k$ be the one that has the type of $k$. It shares only $v$ with
the earlier paths, and so ends at a new vertex of $R$.
\end{proof}

\begin{example}[Yablo, odd Yablo, and a version with even cycles]
\label{n3.4}
In Yablo's graph take $R=v_0\to v_1\to\cdots$. The edge $v_i\to v_j$
has type $1+j$, so the edges from $v_i$ form a mixed fan: every
vertex parity-dominates $R$.

In the odd Yablo graph of Example \ref{n1.1} a walk from $v_i$ to
$v_j$ has the parity of $j-i$, and along a ray the level of a vertex
and its index differ by a constant modulo $2$. So all paths from a
vertex to a ray have the same type. No vertex parity-dominates any
ray, although every vertex dominates every ray.

Now add to the odd Yablo graph, for every $i$, vertices
$x_i,y_i$ and the edges $v_i\to x_i$, $v_i\to y_i$, $x_i\to y_i$ and
$y_i\to x_i$.
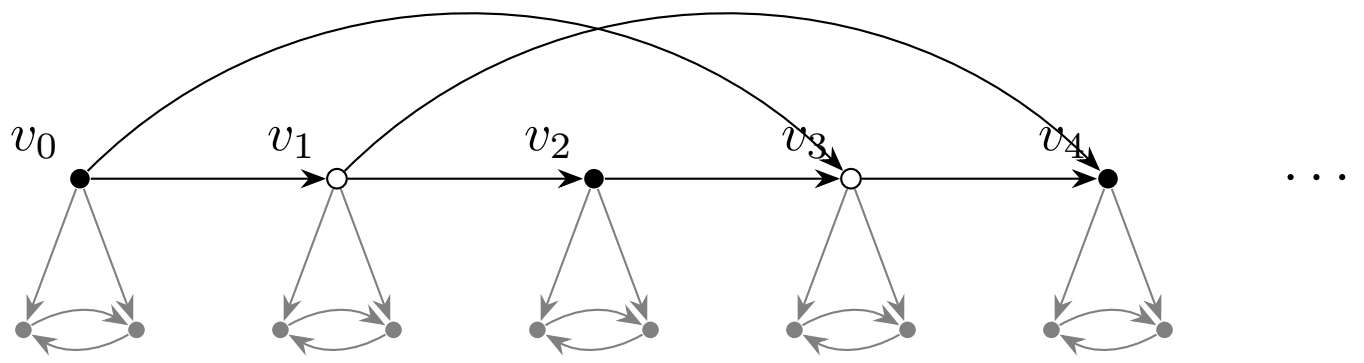
\begin{center}
\begin{tikzpicture}[x=1.7cm,
 sm/.style={circle,fill=gray,inner sep=1.1pt}]
 \foreach \i in {0,2,4} \node[dot] (v\i) at (\i,0){};
 \foreach \i in {1,3} \node[odot] (v\i) at (\i,0){};
 \foreach \i in {0,...,4} \node[above left] at (\i,0){\small$v_{\i}$};
 \node at (4.83,0){$\cdots$};
 \foreach \i/\j in {0/1,1/2,2/3,3/4} \draw[->] (v\i)--(v\j);
 \foreach \i/\j in {0/3,1/4} \draw[->] (v\i) to[bend left=45] (v\j);
 \foreach \i in {0,...,4} {
  \node[sm] (x\i) at (\i-.22,-1){};
  \node[sm] (y\i) at (\i+.22,-1){};
  \draw[->,gray] (v\i)--(x\i); \draw[->,gray] (v\i)--(y\i);
  \draw[->,gray] (x\i) to[bend left=30] (y\i);
  \draw[->,gray] (y\i) to[bend left=30] (x\i); }
\end{tikzpicture}
\end{center}
The digraph has the even cycles $x_i\to y_i\to x_i$ and no odd cycle.
It is not bipartite, because $v_i$ reaches $y_i$ in one step and in
two. The new vertices lie on no ray, and a path between two of the
$v_i$ cannot use them. So fans to rays and
their types are as in the odd Yablo graph: the digraph is parity-safe
and not safe. It has one end. By Theorem \ref{n1.2} it is
kernel-perfect.
\end{example}

\begin{definition}
\label{n3.6}
A ray $R$ \emph{feeds} a ray $R'$, written $R\feeds R'$, if for
every finite set $F$ there is a path in $\Gamma-F$ from a vertex of $R$
to a vertex of $R'$. Write $R\approx R'$ if each feeds the other.
An \emph{end} of $\Gamma$ is an $\approx$-class of rays of $\Gamma$.
\end{definition}

\noo{
\begin{center}
\begin{tikzpicture}
 \draw[->,thick] (0,1.2)--(8,1.2) node[right]{$R$};
 \draw[->,thick] (0,0)--(8,0) node[right]{$R'$};
 \foreach \x in {1,2.6,4.2,5.8} \draw[->,gray] (\x,1.2)--(\x+.7,0);
 \node at (7.2,.6){$\cdots$};
 \node[gate,label=left:{\small$F$}] at (1.4,.6){};
 \node[gate] at (2.9,.65){};
\end{tikzpicture}

\smallskip
{\small $R$ feeds $R'$: no finite set $F$ blocks all connections.}
\end{center}
}

Feeding is always measured in the ambient digraph $\Gamma$, also for rays
that lie in some region of $\Gamma$.

\begin{lemma}[Feeding]
\label{n3.7}
\begin{enumerate}
\item If $R\feeds R'$, then for every finite $F$ and all tails $T$
of $R$ and $T'$ of $R'$ there is a path in $\Gamma-F$ from $T$ to $T'$.
\item $R\feeds R'$ if and only if there are infinitely many
pairwise disjoint paths from $R$ to $R'$.
\item Feeding is reflexive and transitive. A ray and each of its
tails feed each other. So $\approx$ is an equivalence relation, and
a ray and its tails belong to the same end.
\end{enumerate}
\end{lemma}
\begin{proof}
(1) Enlarge $F$ by the finitely many vertices of $R$ outside $T$ and
of $R'$ outside $T'$. (2) A finite set misses one of infinitely many
disjoint paths. Conversely, given finitely many disjoint paths from
$R$ to $R'$, a path from $R$ to $R'$ avoiding all their vertices is
disjoint from them. (3) Paths of length zero show $R\feeds R$, and
the same for a ray and a tail. Let $R\feeds R'\feeds R''$ and let
$F$ be finite. Let $T'$ be a tail of $R'$ that avoids $F$. Take a
path in $\Gamma-F$ from $R$ to a vertex $b$ of $T'$, and by (1) a path in
$\Gamma-F$ from a vertex $b'$ of $R'$ beyond $b$ to $R''$. Together with
the segment of $R'$ from $b$ to $b'$ they form a walk in $\Gamma-F$ from
$R$ to $R''$, which contains a path.
\end{proof}
\paragraph{Remark.} By (2) these are the ends in the sense of \cite{b12}. Ends
in \cite{b10} are coarser, so our results, restricting the number of the former, apply also
to the latter. More specifically, \cite{b10} defines $P\preceq R$, for rays $P,R$, if every
vertex of $P$ reaches some vertex of $R$, and the \emph{end is determined by} $R$ as the
subgraph induced by all vertices that reach $R$. Two rays determine the same end if and
only if $P\preceq R$ and $R\preceq P$; written $P\simeq R$. If $Q$ feeds $R$ then
$P\preceq R$, so $\approx$ is finer than $\simeq$, and the number of ends in the sense of
Definition \ref{n3.6} is not smaller than in the sense of \cite{b10}. 

The two kinds of ends can differ.
Let ray $P=p_0\to p_1\to\cdots$ have  edges $p_i\to r_0$ for all $i$ to a disjoint ray
 $R=r_0\to r_1\to\cdots$. 
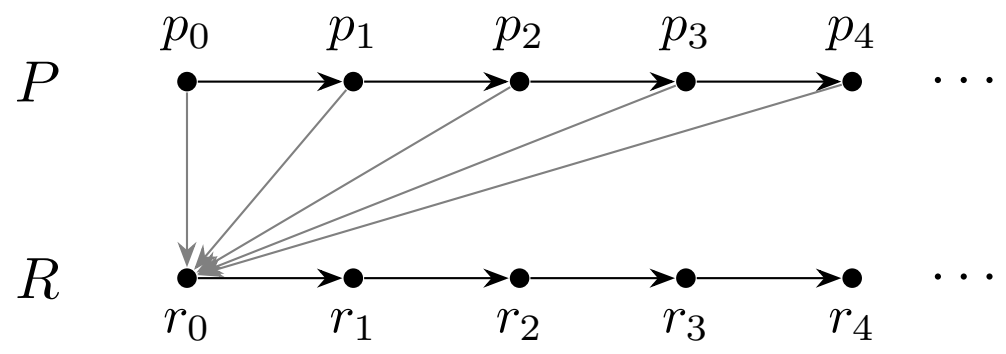
\begin{center}
\begin{tikzpicture}[x=1.1cm]
 \foreach \i in {0,...,4} {
  \node[dot,label=above:{\small$p_{\i}$}] (p\i) at (\i,1.3){};
  \node[dot,label=below:{\small$r_{\i}$}] (r\i) at (\i,0){};}
 \foreach \i/\j in {0/1,1/2,2/3,3/4} {\draw[->] (p\i)--(p\j); \draw[->] (r\i)--(r\j);}
 \foreach \i in {0,...,4} \draw[->,gray] (p\i)--(r0);
 \node at (4.7,1.3){$\cdots$}; \node at (4.7,0){$\cdots$};
 \node at (-.9,1.3){$P$}; \node at (-.9,0){$R$};
\end{tikzpicture}
\end{center}
Then $P\preceq R$, but $P$ does not feed $R$, since the single vertex $r_0$ blocks all
paths. Here $R\not\preceq P$, so the rays determine different ends in both senses. If one
adds the edges $r_i\to p_i$ for even $i$, then $R$ feeds $P$ but still $P\not\approx R$,
while $P\simeq R$.  Both rays lie now in one strong component. In the absence of odd
cycles (e.g., excluding the edges $p_i\to r_0$ for odd $i$), strong components are
kernel-perfect, and this is how we address such situations.

\section{Finite-port regions and the cover lemma}
\label{s4}
This section establishes the central Lemma \ref{n1.3}, allowing to glue solutions of kernel-perefect regions.

\begin{definition}
\label{n4.1}
A set $U\subseteq V(\Gamma)$ is a \emph{finite-port region} of $\Gamma$ if its
set of \emph{ports} $\bd U=\post{\Gamma}(U)\setminus U$ is finite.
\end{definition}

The ports are the vertices outside $U$ that $U$ points to. (They are always taken in the main digraph $\Gamma$.)  Nothing
is assumed about edges \emph{into} $U$. Every successor-closed set
is a finite-port region without ports, and a vertex with finitely
many successors is a finite-port region by itself.
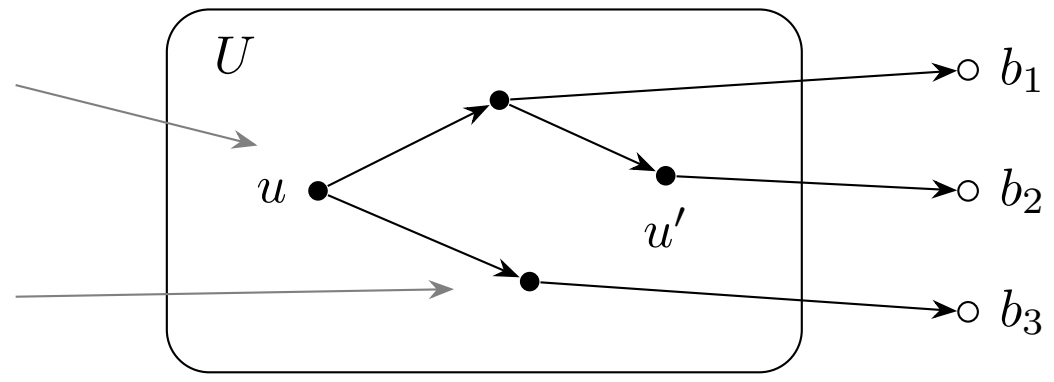
\begin{center}
\begin{tikzpicture}
 \draw[rounded corners=8pt] (0,0) rectangle (4.2,2.4);
 \node at (.45,2.1){\small$U$};
 \node[dot,label=left:{\small$u$}] (u) at (1,1.2){};
 \node[dot] (a) at (2.2,1.8){};
 \node[dot] (b) at (2.4,.6){};
 \node[dot,label=below:{\small$u'$}] (c) at (3.3,1.3){};
 \node[odot,label=right:{\small$b_1$}] (p1) at (5.3,2){};
 \node[odot,label=right:{\small$b_2$}] (p2) at (5.3,1.2){};
 \node[odot,label=right:{\small$b_3$}] (p3) at (5.3,.4){};
 \draw[->] (u)--(a); \draw[->] (u)--(b); \draw[->] (a)--(c);
 \draw[->] (a)--(p1); \draw[->] (b)--(p3); \draw[->] (c)--(p2);
 \draw[->,gray] (-1,1.9)--(.6,1.5); \draw[->,gray] (-1,.5)--(1.9,.55);
\end{tikzpicture}

\smallskip
{\small A region with three ports. The vertex $u$ can reach all of
them; $u'$ can reach only $b_2$. Incoming edges (grey) are
unrestricted.}
\end{center}

\begin{lemma}
\label{n4.2}
\begin{enumerate}
\item The union of two finite-port regions is a finite-port region.
\item If $U$ is a finite-port region of $\Gamma$ and $H$ is an induced
subgraph of $\Gamma$, then $U\cap V(H)$ is a finite-port region of $H$.
\item Let $Z$ be a finite-port region of $\Gamma$, let $v\in Z$, and let
$F\not\ni v$ be finite. The set of vertices reachable from $v$ in
$\Gamma[Z]-F$ is a finite-port region of $\Gamma$, with ports in $\bd Z\cup
F$. We call such a set a \emph{cone} at $v$.
\end{enumerate}
\end{lemma}
\begin{proof}
(1) $\bd(U\cup U')\subseteq\bd U\cup\bd U'$. (2) The ports of
$U\cap V(H)$ in $H$ are among those of $U$ in $\Gamma$. (3) An edge
leaving the cone ends outside $Z$, hence in $\bd Z$, or at a vertex
of $Z$ that is not reachable, which must lie in $F$.
\end{proof}
Viewing a finite-port region $U$ as a device with finitely many ports, once
the truth values of its ports are fixed, the equations inside $U$ no
longer refer to the outside. If $G[U]$ is kernel-perfect they can be
solved for \emph{every} pattern of port values: delete the vertices
that point to a true port, since they must be false, and take a
kernel of what is left. The difficulty is that there may be several
solutions, and the choices made in two regions that look at each
other must fit together.

\begin{example}[Choices that do not fit: an acyclic example]
\label{n4.5}
Take two rays $a_0\to a_1\to\cdots$ and $b_0\to b_1\to\cdots$ and
two more vertices $s\to b_0$ and $t\to a_0$.
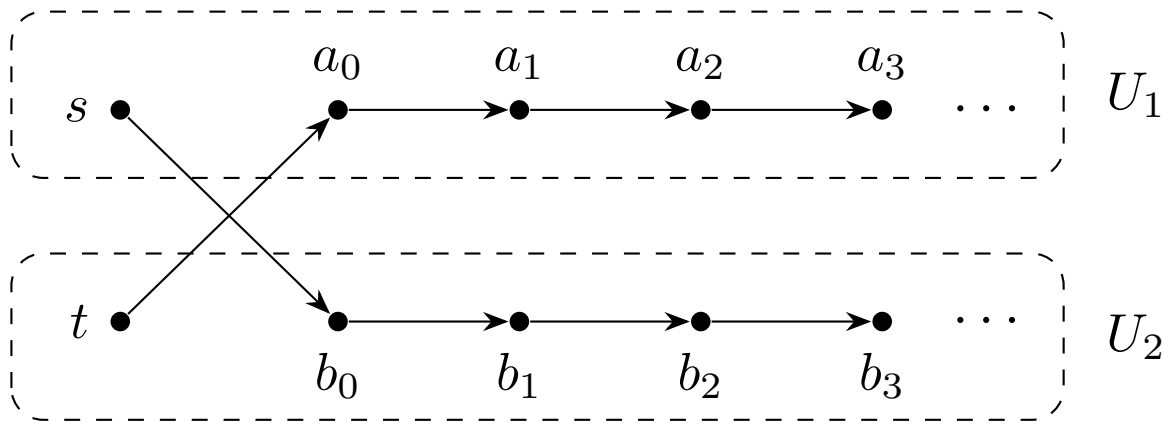
\begin{center}
\begin{tikzpicture}[x=1.2cm]
 \foreach \i in {0,1,2,3} {
  \node[dot,label=above:{\small$a_{\i}$}] (a\i) at (\i+1.2,1.4){};
  \node[dot,label=below:{\small$b_{\i}$}] (b\i) at (\i+1.2,0){}; }
 \foreach \i/\j in {0/1,1/2,2/3} {\draw[->] (a\i)--(a\j); \draw[->] (b\i)--(b\j);}
 \node at (4.8,1.4){$\cdots$}; \node at (4.8,0){$\cdots$};
 \node[dot,label=left:{\small$s$}] (s) at (0,1.4){};
 \node[dot,label=left:{\small$t$}] (t) at (0,0){};
 \draw[->] (s)--(b0); \draw[->] (t)--(a0);
 \draw[dashed,rounded corners=6pt] (-.6,.95) rectangle (5.2,2.05);
 \draw[dashed,rounded corners=6pt] (-.6,-.65) rectangle (5.2,.45);
 \node at (5.6,1.5){\small$U_1$}; \node at (5.6,-.1){\small$U_2$};
\end{tikzpicture}
\end{center}
The region $U_1$ consists of $s$ and the $a$-ray; its only port is
$b_0$. The region $U_2$ consists of $t$ and the $b$-ray; its only
port is $a_0$. A ray has two kernels. So $U_1$ is free to answer its
port by the rule ``$a_0$ is true exactly if $b_0$ is false'', and
$U_2$ by the rule ``$b_0$ is true exactly if $a_0$ is true''.
Together the rules are contradictory, although the digraph has a
kernel. The rule of $U_1$ makes $a_0$ depend on a port that $a_0$
cannot reach. We shall only allow a vertex to depend on the ports it
can reach.
\end{example}

\begin{example}[An example with a cycle]
\label{n4.6}
Let $f\to z$, $z\to y$, $y\to z$, $y\to a$, $a\to f$ and $a\to t$.
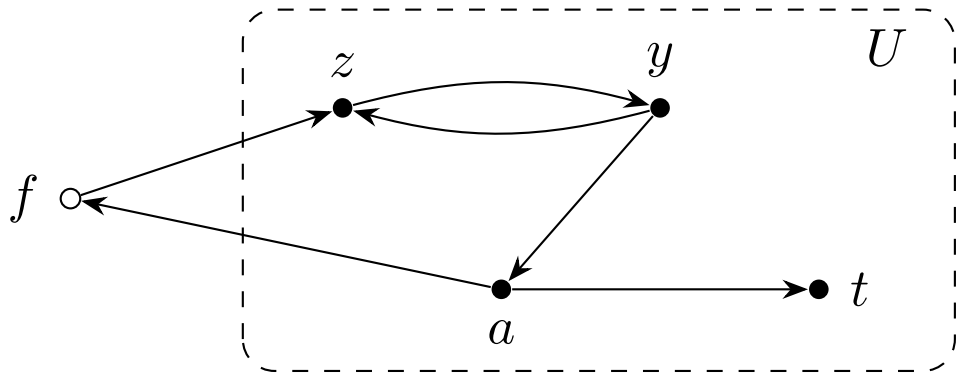
\begin{center}
\begin{tikzpicture}
 \node[odot,label=left:{\small$f$}] (f) at (0,.02){};
 \node[dot,label=above:{\small$z$}] (z) at (1.8,.62){};
 \node[dot,label=above:{\small$y$}] (y) at (3.9,.62){};
 \node[dot,label=below:{\small$a$}] (a) at (2.85,-.58){};
 \node[dot,label=right:{\small$t$}] (t) at (4.95,-.58){};
 \draw[->] (f)--(z);
 \draw[->] (z) to[bend left=15] (y);
 \draw[->] (y) to[bend left=15] (z);
 \draw[->] (y)--(a); \draw[->] (a)--(f); \draw[->] (a)--(t);
 \draw[dashed,rounded corners=6pt] (1.14,-1.12) rectangle (5.85,1.27);
 \node at (5.4,1.02){\small$U$};
\end{tikzpicture}
\end{center}
All cycles are even. Let $U=\{z,y,a,t\}$, with the single port $f$.
The sink $t$ is true, so $a$ is false, whatever $f$ is. Then $z$ and
$y$ form a $2$-cycle on their own, with two solutions. Both reach
$f$. So $U$ may answer by the rule ``$z$ is true exactly if $f$ is
true''. But the equation at $f$ says that $f$ is true exactly if $z$
is false. Reachability no longer excludes the circular dependence,
because $f$ and $z$ lie on a common cycle.

What is wrong with the rule is its direction. The vertex $z$ is an
odd number of steps away from $f$, on the way round the cycle.
Making $f$ true should, if anything, push $z$ towards false. Rules
that respect this are monotone in an order adapted to the
two-colouring of the strong component, and monotone circular systems
have solutions.
\end{example}
In the rest of this section, $G$ is a noc graph. The mentioned order on the assignments is defined as follows. 
Let $S$ be a strong component with its colouring $c$. For $v\in S$
let $\le_v$ be the usual order $0<1$ on $\two$ if $c(v)=0$ and the
reversed order $1<0$ if $c(v)=1$. For assignments $a,a'$ defined on
a set $Y$, write
\[
 a\sqsubseteq_S a'\text{ on }Y\quad\text{if}\quad a(v)\le_v a'(v)\text{ for }v\in Y\cap S\ \text{ and }\ a(v)=a'(v)\text{ for }v\in Y\setminus S.
\]
Going up in this order makes the vertices of colour $0$ truer and
those of colour $1$ falser.

\begin{lemma}[Sign rule]
\label{n4.7}
Let $S$ be a strong component, $v\in S$, $Y\subseteq \post{G}(v)$ be a set of successors of $v$, and
$a\sqsubseteq_S a'$ on $Y$. Put $e(a)=1$ if $a(w)=0$ for all $w\in
Y$, and $e(a)=0$ otherwise. Then $e(a)\le_v e(a')$.
\end{lemma}
\begin{proof}
On the vertices of $Y$ outside $S$ the two assignments agree. The
vertices of $Y$ in $S$ have the colour opposite to that of $v$. If
$c(v)=0$, they are true under $a'$ only if they are true under $a$,
so $e(a)\le e(a')$. If $c(v)=1$, they are true under $a$ only if
they are true under $a'$, so $e(a)\ge e(a')$, which is $e(a)\le_v
e(a')$.
\end{proof}

\paragraph{Solutions and solvers.}
Let $U$ be a finite-port region of $G$. Reaching and strong components always refer to $G$.
A \emph{port pattern} is a map
$\beta:\bd U\to\two$. A \emph{solution of $U$ for $\beta$} is a map
$\alpha:U\to\two$ such that $\alpha\cup\beta$ satisfies the kernel
equation at every vertex of $U$. A \emph{solver} for $U$ is a family
$\sigma=(\sigma_\beta)$ of one solution for each port pattern. For
$u\in U$ let $\bd U(u)$ be the set of ports that $u$ reaches. The
solver is
\begin{itemize}
\item \emph{causal} if $\sigma_\beta(u)$ depends only on
$\beta|\bd U(u)$;
\item \emph{monotone} if for every strong component $S$ of $G$:
$\beta\sqsubseteq_S\beta'$ on $\bd U$ implies
$\sigma_\beta\sqsubseteq_S\sigma_{\beta'}$ on $U\cap S$.
\end{itemize}
Call a solver \emph{good} if it is causal and monotone. In an
acyclic digraph strong components are single vertices and every solver
is monotone. 

The proof of Lemma \ref{n1.3} has three steps: kernel-perfect
regions have good solvers; a good solver survives the addition of
another kernel-perfect region; and a chain of solvers obtained in
this way has a limit.

\begin{lemma}[Layered solving]
\label{n4.8}
Let $W\subseteq V(G)$ induce a kernel-perfect digraph and let
$X\subseteq V(G)\setminus W$ be finite. For $w\in W$ let $X(w)$ be
the set of vertices of $X$ that $w$ reaches. Suppose that for every
$w\in W$ a function $g_w:\two^X\to\two$ is given such that
$g_w(\gamma)$ depends only on $\gamma|X(w)$. Then one can choose,
for every $\gamma\in\two^X$, a map $\omega_\gamma:W\to\two$ such
that for all $w\in W$
\begin{enumerate}
\item $\omega_\gamma(w)=1$ if and only if $g_w(\gamma)=1$ and
$\omega_\gamma(z)=0$ for all $z\in\post{G[W]}(w)$;
\item $\omega_\gamma(w)$ depends only on $\gamma|X(w)$;
\item for every strong component $S$ of $G$: if
$\gamma\sqsubseteq_S\gamma'$ on $X$ implies $g_w(\gamma)\le_w
g_w(\gamma')$ for all $w\in W\cap S$, then
$\gamma\sqsubseteq_S\gamma'$ on $X$ implies
$\omega_\gamma\sqsubseteq_S\omega_{\gamma'}$ on $W\cap S$.
\end{enumerate}
\end{lemma}

Read $g_w(\gamma)=0$ as a veto: under the parameter values
$\gamma$, something outside $W$ forces $w$ to be false.

\begin{proof}
For $\Phi\subseteq X$ let $W_\Phi=\{w\in W:X(w)=\Phi\}$, the
\emph{layer} of $\Phi$. If $w\to z$ is an edge inside $W$ then
$X(z)\subseteq X(w)$: an edge stays in a layer or leads to a layer
with a strictly smaller set. There are finitely many layers.

Call a strong component $S$ \emph{marked} if it meets both $X$ and
$W$. Let $S$ be marked and $x\in X\cap S$. All vertices of $W\cap S$
reach the same vertices, so $W\cap S$ lies in one layer $W_\Phi$,
and $x\in\Phi$. If $w\in W\cap S$ and $w\to z$ with $z\in W_\Phi$,
then $z$ reaches $x$, so $z\in S$: within its layer, $W\cap S$ is
successor-closed. If $S'$ is another marked component meeting
$W_\Phi$, then $S'$ reaches $x\in S$ and symmetrically $S$ reaches
$S'$, so $S'=S$. Hence each layer contains the trace $W\cap S$ of at
most one marked component.

We define the maps $\omega_\gamma$ on the layers $W_\Phi$ by
induction on $|\Phi|$. Fix a pattern $\delta$ on $\Phi$ and consider
all $\gamma$ with $\gamma|\Phi=\delta$ at once; for $w\in W_\Phi$
the value $g_w(\gamma)$ and the values at successors of $w$ in lower
layers are determined by $\delta$. Call $w\in W_\Phi$
\emph{vetoed} if $g_w(\gamma)=0$ or some successor of $w$ in a lower
layer has value $1$.
\begin{itemize}
\item If the layer contains the trace $W\cap S$ of a marked
component, let $Q$ be the set of its vertices that are not vetoed,
and let $\omega_\gamma$ be $1$ on the canonical kernel $\kappa(Q)$
of the bipartite digraph $G[W\cap S]$, with $B_i=\{w:c(w)=i\}$, and
$0$ on the rest of $W\cap S$.
\item On the rest of the layer, veto in addition the vertices with a
successor of value $1$ in $W\cap S$, let $K$ be a kernel of the
subgraph induced by the vertices that are not vetoed, which exists
because $G[W]$ is kernel-perfect, and let $\omega_\gamma$ be $1$ on
$K$ and $0$ elsewhere.
\end{itemize}
All choices depend only on $\delta$, which gives (2). For (1): a
vetoed vertex has value $0$, as it should. A vertex that is not
vetoed has no successor of value $1$ outside its own part of the
layer and satisfies $g_w(\gamma)=1$; within its part, the kernel
property says that it has value $1$ exactly if its successors there
have value $0$.

(3) Let $\gamma\sqsubseteq_S\gamma'$ on $X$. If $S$ is not marked,
either $W\cap S=\varnothing$ or $\gamma=\gamma'$, and there is
nothing to show. Let $S$ be marked, with trace in $W_\Phi$. A
successor of $w\in W\cap S$ in a lower layer reaches no vertex of
$X\cap S$, or it would belong to $S$. By (2) its value is the same
for $\gamma$ and $\gamma'$. So the vetoes on $W\cap S$ differ only
through $g_w$, and by assumption $Q_\gamma\cap B_0\subseteq
Q_{\gamma'}\cap B_0$ and $Q_\gamma\cap B_1\supseteq Q_{\gamma'}\cap
B_1$. Lemma \ref{n2.6}(2) gives
$\omega_\gamma\sqsubseteq_S\omega_{\gamma'}$ on $W\cap S$.
\end{proof}

\begin{lemma}[Kernel-perfect regions have good solvers]
\label{n4.9}
Every finite-port region that induces a kernel-perfect digraph has a
good solver.
\end{lemma}
\begin{proof}
Apply Lemma \ref{n4.8} with $W=U$, $X=\bd U$, and $g_u(\beta)=1$ if
and only if $\beta(b)=0$ for all ports $b\in\post{G}(u)$. These ports
are reached by $u$. Property (1) is the kernel equation at $u$ for
$\omega_\beta\cup\beta$, and (2) is causality. The hypothesis of (3)
holds by the sign rule, applied to the set $Y$ of ports among the
successors of $u$, so (3) is monotonicity.
\end{proof}

Now let $U\subseteq H$ be finite-port regions. Every port of $U$
lies in $H$ or is a port of $H$. A solver $\tau$ for $H$
\emph{extends} a solver $\sigma$ for $U$ if, for every port pattern
$\beta$ of $H$,
\[
 \tau_\beta|U=\sigma_{\beta'}\qquad\text{where }\beta'=(\tau_\beta\cup\beta)|\bd U .
\]
Extension is transitive, as one sees by writing out the definition
twice.

\begin{lemma}[Adding a region]
\label{n4.10}
Let $U$ be a finite-port region with a good solver $\sigma$, and let
$V$ be a finite-port region that induces a kernel-perfect digraph.
Then $H=U\cup V$ has a good solver that extends $\sigma$.
\end{lemma}

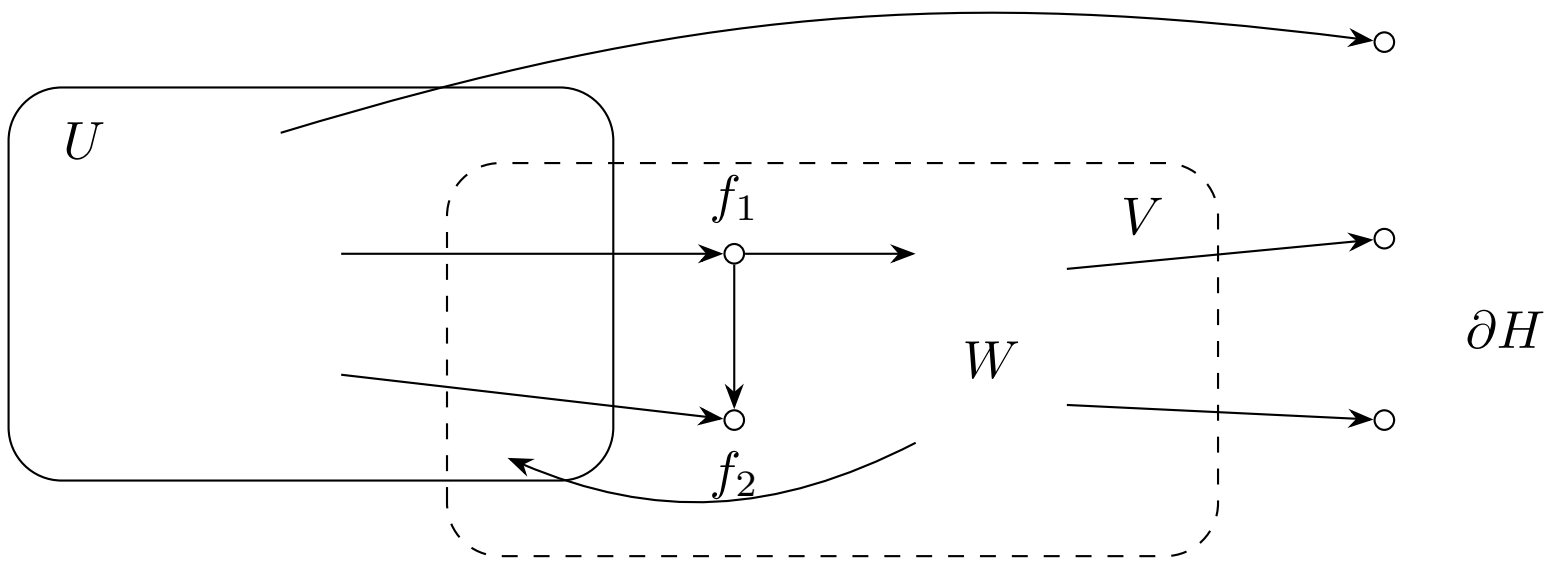
\begin{center}
\begin{tikzpicture}
 \draw[rounded corners=10pt] (0,0) rectangle (4,2.6);
 \node at (.5,2.25){\small$U$};
 \draw[rounded corners=10pt,dashed] (2.9,-.5) rectangle (8,2.1);
 \node at (7.5,1.75){\small$V$};
 \node[odot,label=above:{\small$f_1$}] (f1) at (4.8,1.5){};
 \node[odot,label=below:{\small$f_2$}] (f2) at (4.8,.4){};
 \node at (6.5,.8){\small$W$};
 \node[odot] (x1) at (9.1,1.6){};
 \node[odot] (x2) at (9.1,.4){};
 \node at (9.9,1){\small$\bd H$};
 \node[odot] (x3) at (9.1,2.9){};
 \draw[->] (2.2,1.5)--(f1); \draw[->] (2.2,.7)--(f2);
 \draw[->] (1.8,2.3) to[bend left=12] (x3);
 \draw[->] (f1)--(6,1.5); \draw[->] (7,1.4)--(x1); \draw[->] (7,.5)--(x2);
 \draw[->] (6,.25) to[bend left=25] (3.3,.15);
 \draw[->] (f1)--(f2);
\end{tikzpicture}

\smallskip
{\small $H=U\cup V$. The ports of $U$ are the vertices of
$F=\{f_1,f_2\}$ inside $H$ and some ports of $H$. The rest of $V$ is
$W$.}
\end{center}

\begin{proof}
$H$ is a finite-port region. Put
\[
 F=\bd U\cap H,\qquad X=F\cup\bd H,\qquad W=H\setminus(U\cup F).
\]
Then $\bd U\subseteq X$, the set $X$ is finite, and $W\subseteq V$
induces a kernel-perfect digraph. For a vertex $v$ let $X(v)$ be the
set of vertices of $X$ that $v$ reaches by a path of positive
length. (In Lemma \ref{n4.8}, $X(w)$ admitted empty paths. The two coincide on $W$, since $W\cap X=\emptyset$, but for $f\in F$, now $f\in X(f)$ only if $f$ lies on a cycle.)

\emph{Step 1: pretend the values on $X$ are known.} Let
$\gamma\in\two^X$. On $U$ use the solution $\sigma_{\gamma|\bd U}$.
For $w\in W$ let $g_w(\gamma)=1$ if and only if all successors of
$w$ in $U\cup X$ have value $0$. Each such successor is in $X(w)$,
or is a vertex $u\in U$ whose ports are reached by $w$; as $\sigma$
is causal, $g_w(\gamma)$ depends only on $\gamma|X(w)$, although $w$
may have infinitely many successors in $U$. Lemma \ref{n4.8} gives
maps $\omega_\gamma$ on $W$. Let $a_\gamma$ be the assignment on
$H\cup\bd H$ that is $\gamma$ on $X$, $\sigma_{\gamma|\bd U}$ on $U$
and $\omega_\gamma$ on $W$. Then:
\begin{enumerate}
\renewcommand{\labelenumi}{(\alph{enumi})}
\item $a_\gamma$ satisfies the kernel equation at every vertex of
$U$ and of $W$;
\item $a_\gamma(v)$ depends only on $\gamma|X(v)$, for $v\in U\cup W$;
\item for every strong component $S$: if $\gamma\sqsubseteq_S\gamma'$
on $X$, then $a_\gamma\sqsubseteq_S a_{\gamma'}$ on $(H\cup\bd
H)\cap S$.
\end{enumerate}
We also note
\begin{enumerate}
\renewcommand{\labelenumi}{(\alph{enumi})}
\setcounter{enumi}{3}
\item if $\gamma\sqsubseteq_S\gamma'$ on $X$ and $v\notin S$ is
reached from $S$, then $a_\gamma(v)=a_{\gamma'}(v)$.
\end{enumerate}
Indeed such a $v$ reaches no vertex of $S$; apply (b), or
$\gamma(v)=\gamma'(v)$ if $v\in X$. For (c): on $X\cap S$ this is
the assumption, and on $U\cap S$ it holds because $\sigma$ is
monotone. Let $w\in W\cap S$ and let $Y$ be the set of its
successors in $U\cup X$. By what was just said and (d),
$a_\gamma\sqsubseteq_S a_{\gamma'}$ on $Y$. The sign rule gives
$g_w(\gamma)\le_w g_w(\gamma')$, and Lemma \ref{n4.8}(3) gives the
claim on $W\cap S$.

\emph{Step 2: the equations at $F$.} The vertices of $F$ belong to
$H$, so their equations must hold too. All successors of $f\in F$
lie in $H\cup\bd H$. So the equation at $f$ reads
\[
 \gamma(f)=h_f(\gamma),\qquad h_f(\gamma)=1\iff a_\gamma(z)=0\text{ for all }z\in\post{G}(f),
\]
and $h_f(\gamma)$ depends only on $\gamma|X(f)$: for successors of
$f$ in $U\cup W$ by (b), and successors of $f$ in $X$ belong to
$X(f)$. Group the vertices of $F$ into classes of vertices that
reach each other, and list the classes as $K_1,\dots,K_m$ so that a
class reaches only itself and classes listed before it. Given a port
pattern $\beta$ of $H$, put $\gamma=\beta$ on $\bd H$ and determine
$\gamma$ on $K_1,K_2,\dots$ in turn. When the turn of $K=K_j$ comes,
$\gamma$ is known on $\bd H$ and on the earlier classes, and
$X(f)\subseteq\bd H\cup K_1\cup\dots\cup K_j$ for $f\in K$. For
$\delta\in\two^K$ let $\gamma_\delta$ be the assignment on $X$ that
has these known values, is $\delta$ on $K$, and is $0$ on the later
classes, which do not matter. We need a $\delta$ with
$\Psi(\delta)=\delta$, where
\[
 \Psi:\two^K\to\two^K,\qquad \Psi(\delta)(f)=h_f(\gamma_\delta).
\]

Let $S$ be the strong component containing $K$. If
$\delta\sqsubseteq_S\delta'$ on $K$, then
$\gamma_\delta\sqsubseteq_S\gamma_{\delta'}$ on $X$. By (c) and (d)
the assignments $a_{\gamma_\delta}$ and $a_{\gamma_{\delta'}}$
satisfy $\sqsubseteq_S$ on the set of successors of any $f\in K$,
and the sign rule gives $\Psi(\delta)\sqsubseteq_S\Psi(\delta')$. So
$\Psi$ is a monotone map of the finite partially ordered set
$(\two^K,\sqsubseteq_S)$, which has a largest element $\delta_0$.
Then $\delta_0\sqsupseteq\Psi(\delta_0)\sqsupseteq\Psi^2(\delta_0)
\sqsupseteq\cdots$, and this sequence becomes constant at a fixed
point $\delta^*$ of $\Psi$. It is the largest one: a fixed point
$\delta$ satisfies $\delta=\Psi^n(\delta)\sqsubseteq_S
\Psi^n(\delta_0)$. Let $\gamma|K=\delta^*$.

(If $K$ is a single vertex on no cycle, then $f\notin X(f)$, the map
$\Psi$ is constant, and this step just evaluates $h_f$. In an
acyclic digraph this is all that happens.)

Write $\gamma^\beta$ for the resulting $\gamma$ and let
$\tau_\beta=a_{\gamma^\beta}|H$.

\emph{Step 3: checking.} $\tau_\beta$ is a solution of $H$ for
$\beta$, by (a) and Step 2. It extends $\sigma$ by construction.

\emph{Causal.} By induction along the list of classes, $\gamma^\beta$
on $K_j$ depends only on the values of $\beta$ at the ports of $H$
reachable from $K_j$: the map $\Psi$ for $K_j$ is determined by
$\gamma^\beta$ on $\bd H$ and on earlier classes, restricted to what
$K_j$ reaches. With (b), $\tau_\beta(v)$ depends only on the values
of $\beta$ at the ports of $H$ that $v$ reaches.

\emph{Monotone.} Let $S$ be a strong component and
$\beta\sqsubseteq_S\beta'$ on $\bd H$. Write $\gamma=\gamma^\beta$
and $\gamma'=\gamma^{\beta'}$. A class of $F$ that does not reach
$S$ reaches no port in $S$, so by causality $\gamma=\gamma'$ on it.

Suppose first that $K=F\cap S$ is nonempty. It is one of the
classes, and the classes listed before it do not reach $S$. Let
$\Psi,\Psi'$ be the maps of $K$ for $\beta$ and $\beta'$, built from
assignments $\gamma_\delta,\gamma'_\delta$, and let
$\delta^*,\delta'^*$ be their largest fixed points. Then
$\gamma_\delta\sqsubseteq_S\gamma'_\delta$ on $X$, so as in Step 2,
$\Psi(\delta)\sqsubseteq_S\Psi'(\delta)$ for every $\delta$. Hence
$\delta^*=\Psi(\delta^*)\sqsubseteq_S\Psi'(\delta^*)$, and applying
the monotone map $\Psi'$ repeatedly,
$\delta^*\sqsubseteq_S\Psi'^n(\delta^*)\sqsubseteq_S
\Psi'^n(\delta_0)$ for all $n$. So $\delta^*\sqsubseteq_S\delta'^*$.

In either case let $\hat\gamma$ be the assignment on $X$ that agrees
with $\gamma$ on the vertices of $X$ reached from $S$ and with
$\gamma'$ on the others. A vertex of $X$ reached from $S$ lies in
$\bd H$, in $K$, or in a class that does not reach $S$. So
$\hat\gamma\sqsubseteq_S\gamma'$ on $X$, and (c) gives
$a_{\hat\gamma}\sqsubseteq_S a_{\gamma'}$ on $H\cap S$. Moreover
$a_{\hat\gamma}=a_\gamma$ on $H\cap S$: on $F\cap S$ because
$\hat\gamma=\gamma$ there, and on $(U\cup W)\cap S$ by (b), since a
vertex of $S$ reaches only vertices reached from $S$. Hence
$\tau_\beta\sqsubseteq_S\tau_{\beta'}$ on $H\cap S$.
\end{proof}
We have now all the pieces for the aim of this section.
\begin{restatedL}{n1.3}[Finite-port cover]
If every vertex of a noc graph $G$ lies
in a finite-port region that induces a kernel-perfect digraph, then
$G$ is kernel-perfect.
\end{restatedL}
\begin{proof} 
The hypothesis passes to induced subgraphs, by Lemma \ref{n4.2}(2).
By Proposition \ref{n2.3} it therefore suffices to show that a
nonempty noc graph $G$ satisfying the hypothesis has a nonempty
semikernel. We find a nonempty successor-closed set $D$ such that
$G[D]$ has a kernel; this is enough by Lemma \ref{n2.2}(2).

Let $U_0$ be a kernel-perfect finite-port region
containing some vertex, and $\sigma^0$ a good solver for it
(Lemma \ref{n4.9}). Given a finite-port region $U_n$ with a good
solver $\sigma^n$, choose for each of its finitely many ports a
kernel-perfect finite-port region containing that port, and add
these regions to $U_n$ one at a time, using Lemma \ref{n4.10}. The
result is a finite-port region $U_{n+1}\supseteq U_n\cup\bd U_n$
with a good solver $\sigma^{n+1}$ that extends $\sigma^n$. Let
$D=\bigcup_n U_n$. It is successor-closed, because
$\post{G}(U_n)\subseteq U_{n+1}$.
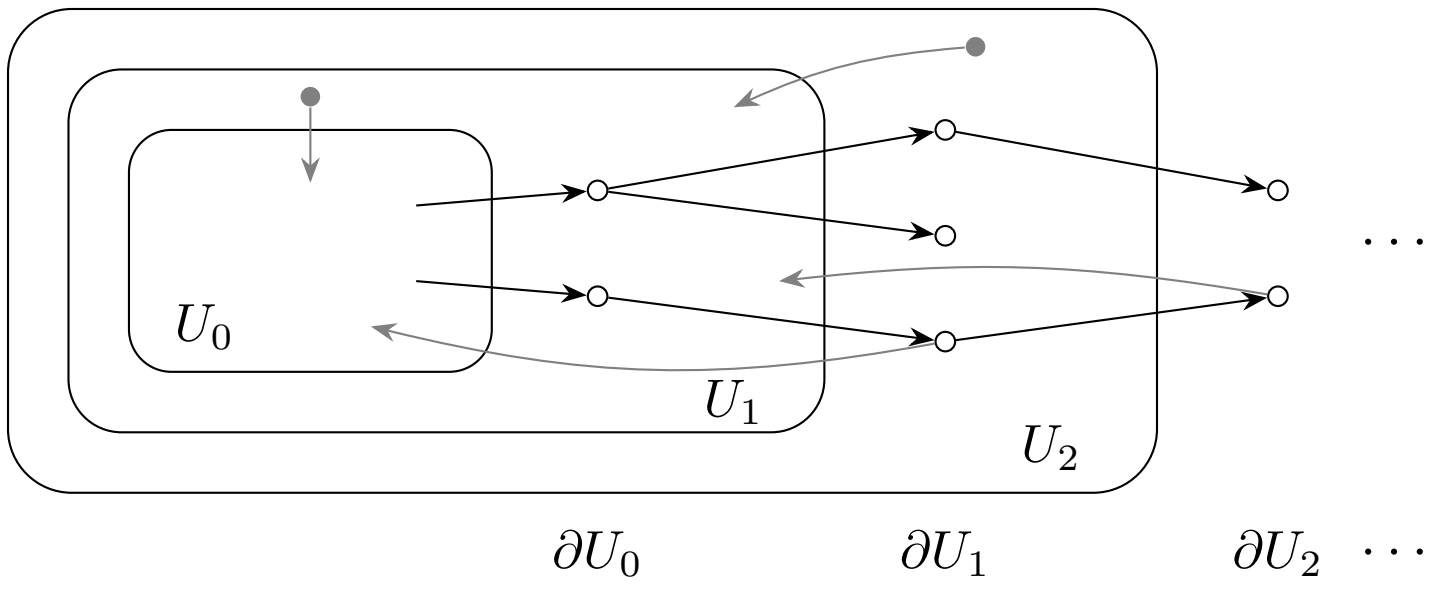
\begin{center}
\begin{tikzpicture}
 \draw[rounded corners=8pt] (0,0) rectangle (2.4,1.6);
 \draw[rounded corners=10pt] (-.4,-.4) rectangle (4.6,2.0);
 \draw[rounded corners=12pt] (-.8,-.8) rectangle (6.8,2.4);
 \node at (.5,.3){\small$U_0$}; \node at (4.0,-.2){\small$U_1$};
 \node at (6.1,-.5){\small$U_2$};
 \node[odot] (p1) at (3.1,1.2){}; \node[odot] (p2) at (3.1,.5){};
 \node[odot] (q1) at (5.4,1.6){}; \node[odot] (q2) at (5.4,.9){};
 \node[odot] (q3) at (5.4,.2){};
 \node[odot] (r1) at (7.6,1.2){}; \node[odot] (r2) at (7.6,.5){};
 \draw[->] (1.9,1.1)--(p1); \draw[->] (1.9,.6)--(p2);
 \draw[->] (p1)--(q1); \draw[->] (p1)--(q2); \draw[->] (p2)--(q3);
 \draw[->] (q1)--(r1); \draw[->] (q3)--(r2);
 \draw[->,gray] (q3) to[bend left=12] (1.6,.3);
 \draw[->,gray] (r2) to[bend right=8] (4.3,.6);
 \node[dot,gray] (u1) at (1.2,1.82){}; \draw[->,gray] (u1)--(1.2,1.25);
 \node[dot,gray] (u3) at (5.6,2.15){}; \draw[->,gray] (u3) to[bend right=10] (4.0,1.75);
 \node at (8.4,.85){$\cdots$};
 \node at (3.1,-1.2){\small$\bd U_0$};
 \node at (5.4,-1.2){\small$\bd U_1$};
 \node at (7.6,-1.2){\small$\bd U_2$};
 \node at (8.4,-1.2){$\cdots$};
\end{tikzpicture}

\smallskip
{\small The ports of $U_n$ lie inside $U_{n+1}$. A state of the ports
of $U_{n+1}$ determines a state of the ports of $U_n$. Edges back
into $U_n$ (grey) are unrestricted.}
\end{center}
Let \emph{states} be elements of finite nonempty $P_n=\two^{\bd U_n}$. Since
$\bd U_n\subseteq U_{n+1}$, every state $\beta\in P_{n+1}$ determines a state
$\pi_n(\beta)=\sigma^{n+1}_\beta|\bd U_n\in P_n$. By K\H{o}nig's lemma there is a sequence
of states $\beta_n\in P_n$ with $\pi_n(\beta_{n+1})=\beta_n$ for all $n$: the finite
sequences $(\beta_0,\dots,\beta_n)$ with $\beta_i=\pi_i(\beta_{i+1})$ form an infinite
tree with finite levels, since every $\beta_n\in P_n$ generates one. Because
$\sigma^{n+1}$ extends $\sigma^n$,
\[
 \sigma^{n+1}_{\beta_{n+1}}|U_n=\sigma^n_{\beta'}\quad\text{with}\quad
 \beta'=\sigma^{n+1}_{\beta_{n+1}}|\bd U_n=\pi_n(\beta_{n+1})=\beta_n .
\]
So the solutions $\sigma^n_{\beta_n}$ agree where they overlap, and their union
$\kappa\in \{0,1\}^D$, with $\kappa|\bd U_n = \beta_n$.  Hence on $U_n\cup\bd U_n$ the
assignment $\kappa$ is $\sigma^n_{\beta_n}\cup\beta_n$, which satisfies the kernel
equation at every vertex of $U_n$. Every vertex of $D$ lies in some $U_n$. So $\kappa$
determines a kernel of $G[D]$.
\end{proof}

\section{Ladder rays: paradox travels along rays}
\label{s5}

\subsection{Collisions}
\label{s5.1}

Write $u\prec v$ if $\Gamma$ contains paths from $u$ to $v$ of both
parities, and call it a \emph{collision}. In a noc graph $G$, by the parity principle, walks of both parities suffice,
and therefore
\begin{equation}
\tag{$*$}\label{e1}
 u\text{ reaches }v,\quad v\prec w,\quad w\text{ reaches }x
 \quad\Longrightarrow\quad u\prec x .
\end{equation}
In a noc graph, if $u\prec v$ then $u$ and $v$ lie in different strong components,
by Lemma \ref{n2.1}(3); in particular $u\ne v$.

\begin{lemma}[Collision]
\label{n5.1}
Let $C$ be a core without odd cycles and $r$ a vertex of $C$. There
is a vertex $z$ of $C$ with $r\prec z$ in $C$.
\end{lemma}

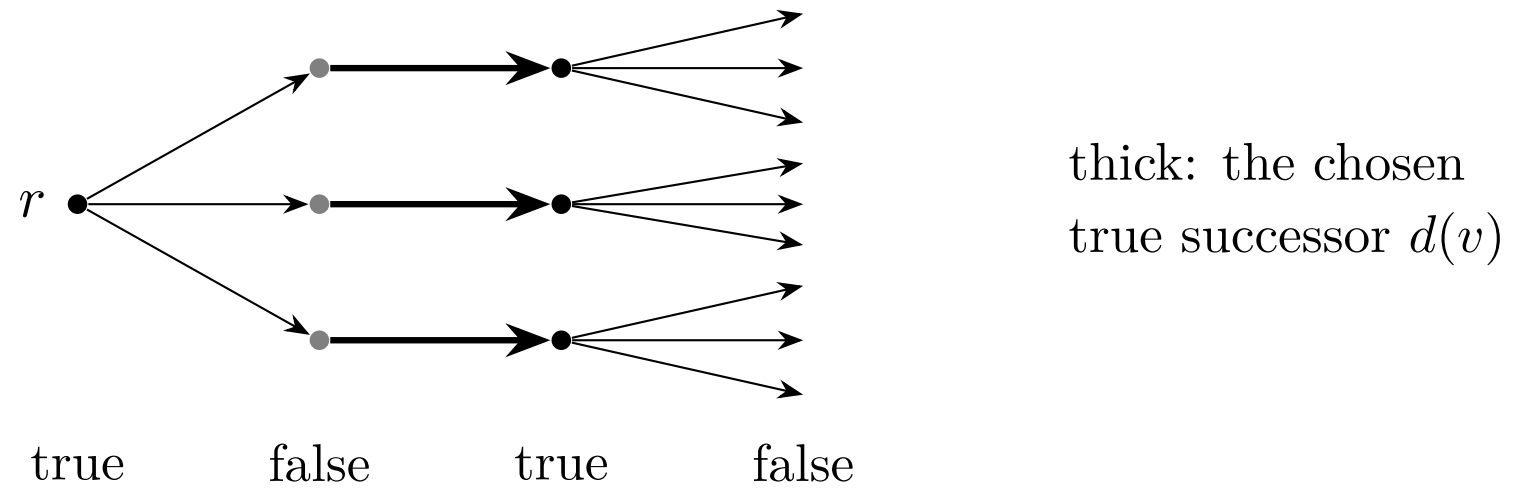
\begin{center}
\begin{tikzpicture}[x=1.6cm,y=.9cm]
 \node[dot,label=left:{$r$}] (r) at (0,0){};
 \foreach \y/\n in {1/a,0/b,-1/c} \node[dot,gray] (\n) at (1,\y){};
 \foreach \y/\n in {1/a2,0/b2,-1/c2} \node[dot] (\n) at (2,\y){};
 \foreach \n in {a,b,c} \draw[->] (r)--(\n);
 \foreach \n in {a,b,c} \draw[->,very thick] (\n)--(\n2);
 \foreach \y in {1.4,1,.6} \draw[->] (a2)--(3,\y);
 \foreach \y in {.3,0,-.3} \draw[->] (b2)--(3,\y);
 \foreach \y in {-.6,-1,-1.4} \draw[->] (c2)--(3,\y);
 \node at (0,-1.9){\small true};
 \node at (1,-1.9){\small false};
 \node at (2,-1.9){\small true};
 \node at (3,-1.9){\small false};
 \node[align=left] at (5,0){\small thick: the chosen\\ \small true successor $d(v)$};
\end{tikzpicture}
\end{center}

\begin{proof}
Try to build a semikernel containing $r$. Declare $r$ true. Every
successor of a true vertex must be false; here there is no choice.
Every false vertex needs a true successor; here there is a choice,
and we make it once and for all: since $C$ has no sinks, pick a
successor $d(v)$ for every vertex $v$. Call a walk
$r=x_0\to x_1\to\cdots\to x_\ell$ \emph{alternating} if
$x_{t+1}=d(x_t)$ for all odd $t$. Let $E$ be the set of last
vertices of alternating walks of even length, the vertices that the
attempt makes true, and $O$ the same for odd length, the vertices it
makes false. An even alternating walk can be prolonged by any edge,
and an odd one by the chosen edge, so
\[
 \post{C}(E)\subseteq O,\qquad d(O)\subseteq E .
\]
If $E$ and $O$ were disjoint, then $E$ would be independent, by the
first inclusion, and every successor of $E$ would have a successor
in $E$, by the second. So $E$ would be a semikernel containing $r$.
A core has none. Hence some $z$ is reached from $r$ by walks of both
parities.
\end{proof}

\subsection{Essential vertices and ladder rays}\label{s5.2}

\begin{definition}\label{n4.11}
A vertex $v$ of $\Gamma$ is \emph{essential} in $\Gamma$ if no finite-port
region of $\Gamma$ that contains $v$ induces a kernel-perfect digraph.
Write $E_\Gamma$ for the set of essential vertices.
\end{definition}

In Yablo's graph every vertex is essential: a finite-port region
around $v_i$ contains all but finitely many of the $v_j$. If Yablo's
graph is extended by a vertex $s$ with the single edge $s\to v_0$,
then $s$ is not essential, since $\{s\}$ is a kernel-perfect
finite-port region. The paradox is ``at'' the $v_i$ and not at $s$.

\begin{lemma}[Essential vertices]
\label{n4.12}
\begin{enumerate}
\item If $H$ is an induced subgraph of $\Gamma$, then $E_H\subseteq E_\Gamma$.
\item A vertex $v$ is essential if and only if no cone at $v$
induces a kernel-perfect digraph. In detail: if $v$ is essential, $Z$
is a finite-port region containing $v$, and $F\not\ni v$ is finite,
then the set of vertices reachable from $v$ in $\Gamma[Z]-F$ induces a
digraph that is not kernel-perfect.
\item In a core $C$ without odd cycles, every vertex reaches an essential vertex of
$C$.
\end{enumerate}
\end{lemma}
\begin{proof}
(1) If $v\in V(H)$ lies in a finite-port region $U$ of $\Gamma$ with
$\Gamma[U]$ kernel-perfect, then $U\cap V(H)$ is such a region of $H$.

(2) Cones are finite-port regions by Lemma \ref{n4.2}(3).
Conversely, if $U\ni v$ is a kernel-perfect finite-port region, the
cone of vertices reachable from $v$ in $\Gamma-\bd U$ is contained in
$U$.

(3) Let $D$ be the future of a vertex in $C$. By
Corollary \ref{n2.4}, $C[D]$ is not kernel-perfect. By
Lemma \ref{n1.3} some $w\in D$ is essential in $C[D]$, and by (1)
it is essential in $C$.
\end{proof}

A kernel-perfect digraph has no essential vertices. By (3) and
Corollary \ref{n2.4}, a noc graph that is not kernel-perfect has some. The rest of this section shows that it has  
a special ladder ray with infinitely many essential vertices.

\begin{definition}
\label{n5.2}
A \emph{ladder ray} of $\Gamma$ is a ray that contains essential vertices
$b_1,b_2,b_3,\dots$ of $\Gamma$, in this order, with $b_i\prec b_j$
whenever $i<j$. The $b_i$ are \emph{checkpoints} of the ladder ray.
\end{definition}

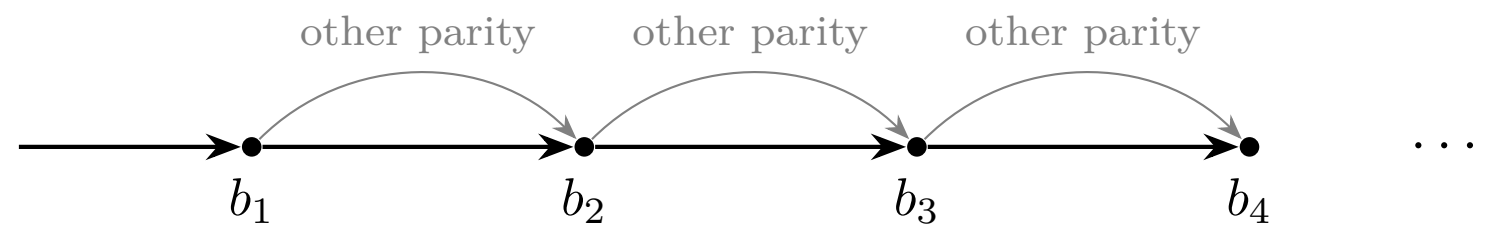
\begin{center}
\begin{tikzpicture}[x=2.2cm]
 \foreach \i in {1,...,4} \node[dot,label=below:{\small$b_{\i}$}] (b\i) at (\i,0){};
 \node at (4.6,0){$\cdots$};
 \draw[->,thick] (.3,0)--(b1);
 \foreach \i/\j in {1/2,2/3,3/4}{
  \draw[->,thick] (b\i)--(b\j);
  \draw[->,gray] (b\i) to[bend left=45] node[above]{\scriptsize other parity}(b\j);}
\end{tikzpicture}

\smallskip
{\small A ladder ray (thick). Each checkpoint is essential and
reaches the next one also by a path of the other parity, which need
not be disjoint from the ray.}
\end{center}

\begin{lemma}
\label{n5.3}
\begin{enumerate}
\item A tail of a ladder ray is a ladder ray, and belongs to the
same end.
\item \emph{Splicing.} In a noc graph, let a ray $R$ pass through vertices
$p_0,q_0,p_1,q_1,\dots$ in this order, where each $p_n$ is
essential and $p_n\prec q_n$. Then $R$ is a ladder ray with
checkpoints $p_n$.
\item In a noc graph, no tail of a ladder ray lies in a single strong component.
\end{enumerate}
\end{lemma}
\begin{proof}
(1) Drop the checkpoints outside the tail. (2) For $n<m$ the vertex
$q_n$ reaches $p_m$ along $R$, so $p_n\prec p_m$ by (\ref{e1}).
(3) Checkpoints $b_i\prec b_j$ lie in different strong components.
\end{proof}
In an acyclic digraph one finds a ladder ray by walking from collision
to collision. With cycles the walk might return to where it has
been. The next lemma takes one step avoiding
any given set of vertices.

\begin{lemma}[Ladder step]
\label{n5.4}
Let $H$ be an induced subgraph of a noc graph $G$, $v\in V(H)$, 
$A\subseteq V(H)\setminus\{v\}$, and suppose that the set
$Y$ of vertices reachable from $v$ in $H-A$ induces a digraph that is
not kernel-perfect. Then there is a path $P$ in $H[Y]$ from $v$ to a
vertex $q$, passing through a vertex $p$ before $q$, such that
\begin{enumerate}
\item $p$ is essential in $G$ and $p\prec q$;
\item the set of vertices reachable from $q$ in
$H-(A\cup V(P)\setminus\{q\})$ induces a digraph that is not
kernel-perfect.
\end{enumerate}
\end{lemma}

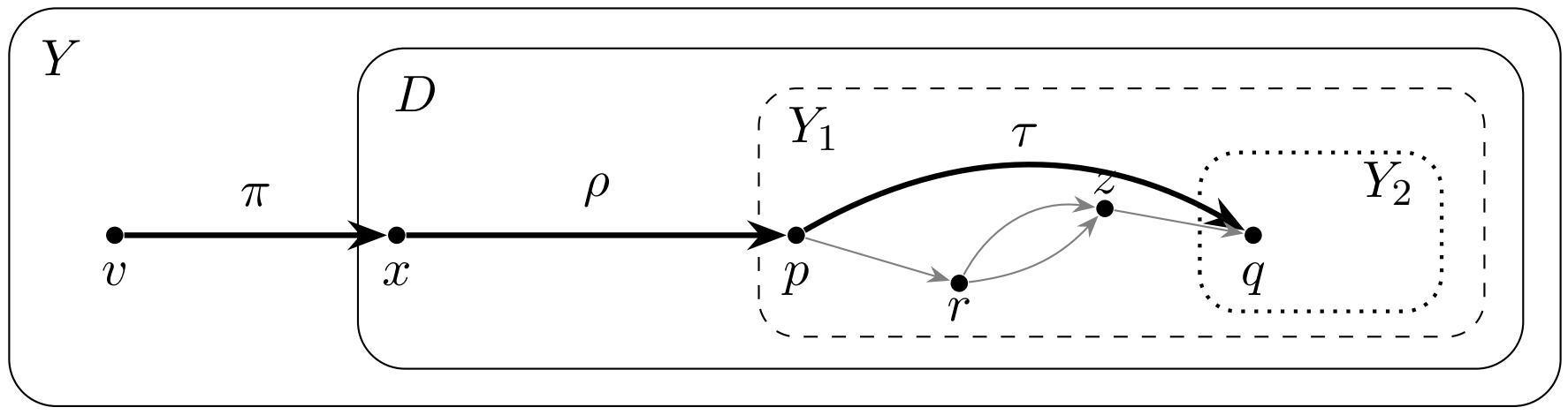
\begin{center}
\begin{tikzpicture}
 \draw[rounded corners=10pt] (-.79,-1.28) rectangle (10.82,1.7);
 \node at (-.4,1.33){$Y$};
 \draw[rounded corners=10pt] (1.82,-1) rectangle (10.54,1.4);
 \node at (2.25,1.06){$D$};
 \draw[rounded corners=8pt,dashed] (4.82,-.76) rectangle (10.25,1.1);
 \node at (5.23,.8){$Y_1$};
 \draw[rounded corners=8pt,dotted,thick] (8.12,-.57) rectangle (9.93,.62);
 \node at (9.53,.39){$Y_2$};
 \node[dot,label=below:{$v$}] (v) at (0,0){};
 \node[dot,label=below:{$x$}] (x) at (2.11,0){};
 \node[dot,label=below:{$p$}] (p) at (5.1,0){};
 \node[dot,label=below:{$q$}] (q) at (8.52,0){};
 \node[dot,label={[inner sep=1pt]below:{$r$}}] (r) at (6.32,-.36){};
 \node[dot,label={[inner sep=1pt]above:{$z$}}] (z) at (7.41,.2){};
 \draw[->,very thick] (v)--node[above=2pt]{$\pi$}(x);
 \draw[->,very thick] (x)--node[above=2pt]{$\rho$}(p);
 \draw[->,very thick] (p) to[bend left=30] node[above]{$\tau$}(q);
 \draw[->,gray] (p)--(r);
 \draw[->,gray] (r) to[bend left=35] (z);
 \draw[->,gray] (r) to[bend right=20] (z);
 \draw[->,gray] (z)--(q);
\end{tikzpicture}

\smallskip
{\small One ladder step. The path $P=\pi\rho\tau$ is thick. Grey: $p$
reaches $r$, which reaches $z$ in both parities, and $z$ reaches
$q$; so $p\prec q$.}
\end{center}

\begin{proof}
$H[Y]$ contains an induced core $C$. Take a path $\pi$ in $H[Y]$
from $v$ to a vertex of $C$ and stop it at its first vertex $x$ in
$C$. Let $D$ be the future of $x$ in $C$; it induces a core, and
$\pi$ meets $D$ only in $x$. By Lemma \ref{n4.12}(3) there are a
vertex $p$ that is essential in $C[D]$ and a path $\rho$ in $D$
from $x$ to $p$. By Lemma \ref{n4.12}(1), $p$ is essential in $G$.

Let $Y_1$ be the set of vertices reachable from $p$ in
$C[D]-(V(\rho)\setminus\{p\})$. It is a cone at the essential vertex
$p$ of $C[D]$, so it contains an induced core $C_1$. Pick $r\in
C_1$. By Lemma \ref{n5.1} there is $z\in C_1$ with $r\prec z$, and
by Lemma \ref{n4.12}(3) $z$ reaches a vertex $q$ that is essential
in $C_1$. As $p$ reaches $r$, we get $p\prec q$ from (\ref{e1}). Let
$\tau$ be a path in $Y_1$ from $p$ to $q$, and $P=\pi\rho\tau$. This
is a path: $\rho$ and $\tau$ lie in $D$, which $\pi$ meets only in
$x$, and $\tau$ lies in $Y_1$, which $\rho$ meets only in $p$.

Finally let $Y_2$ be the set of vertices reachable from $q$ in
$C_1-(V(\tau)\setminus\{q\})$. It is a cone at the essential vertex
$q$ of $C_1$, so it induces a digraph that is not kernel-perfect. It
avoids $A$ and all vertices of $P$ other than $q$, because
$C_1\subseteq Y_1\subseteq D\subseteq Y$. So it is contained in the
set in (2), which is therefore not kernel-perfect either.
\end{proof}

\begin{proposition}
\label{n5.5}
Let $G$ be a noc graph and $Z\subseteq V(G)$. If $G[Z]$ is not kernel-perfect, 
$Z$ contains some ladder ray of $G$.
\end{proposition}
\begin{proof}
Let $H=G[Z]$, let $q_{-1}$ be a vertex of an induced core of $H$,
and $A_0=\varnothing$. The future of $q_{-1}$ in $H$ contains a
core, so Lemma \ref{n5.4} applies to $q_{-1}$ and $A_0$. It gives a
path $P_0$ from $q_{-1}$ through $p_0$ to $q_0$, and its conclusion
(2) is the hypothesis of the lemma for $q_0$ and
$A_1=V(P_0)\setminus\{q_0\}$. Continuing with
$A_{n+1}=A_n\cup(V(P_n)\setminus\{q_n\})$, we get paths $P_n$ from
$q_{n-1}$ through $p_n$ to $q_n$, where $P_n$ avoids $A_n$, that is,
all vertices of the earlier paths except its first vertex. Their
concatenation is a ray in $Z$ through $p_0,q_0,p_1,q_1,\dots$, and a
ladder ray by Lemma \ref{n5.3}(2).
\end{proof}
A kernel-perfect digraph has no essential vertex and hence no ladder ray. Thus 

\begin{proposition}
A noc graph is kernel-perfect iff it has no ladder ray.
\end{proposition}

\begin{definition}
\label{n5.6}
A \emph{ladder end} is an end that contains a ladder ray. It
\emph{lives in} a set $Z$ if it contains a ladder ray that lies in
$Z$.
\end{definition}
Two parity conditions appear in this picture and should be kept
apart. On a ladder ray, each checkpoint reaches every later one by
paths of both parities, but these paths may all run through the same
vertices. A mixed fan asks for infinitely many \emph{disjoint} paths
of both types from one vertex. Parity-safety says that the first
kind of parity conflict, which every paradox has, is not accompanied
infinitely often by the second.

\section{Gates and the diagonal comb}\label{s6}

\subsection{Passing a gate}\label{s6.1}
Here is the central picture. A ray $L$ lies in a region $Z$, and an
essential vertex $q$ does not parity-dominate it. So there is a
finite gate $F$ such that all paths from $q$ to a tail $T$ of $L$
that avoid $F$ have the same type. Let $Z^*$ be what $q$ can reach
in $Z$ without passing the gate. It falls into two parts: the
vertices that can still get to $T$ without passing the gate, and the
vertices that cannot.

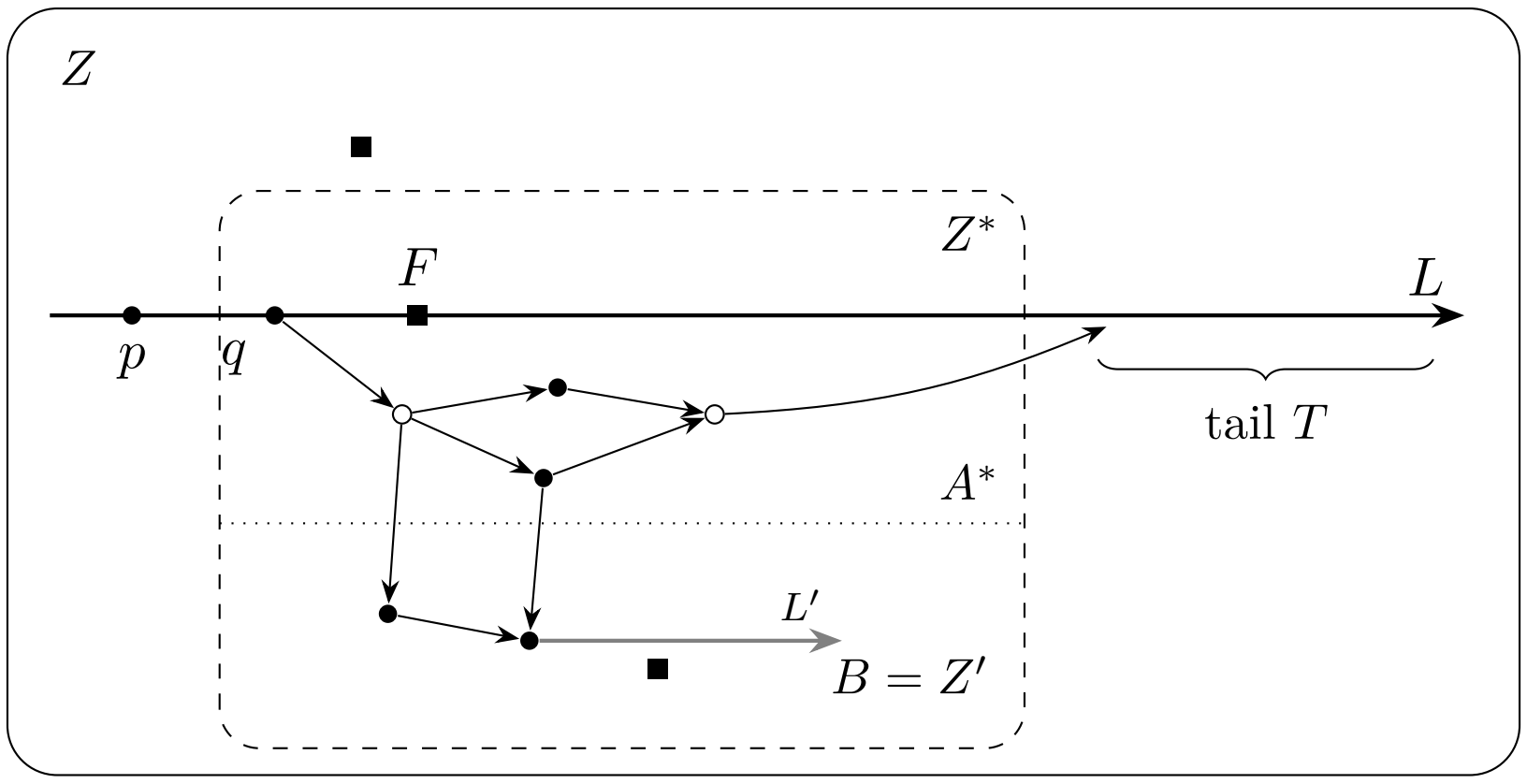
\begin{center}
\begin{tikzpicture}
 \draw[rounded corners=10pt] (-.3,-3.25) rectangle (10.39,2.17);
 \node at (.2,1.75){\small$Z$};
 \draw[->,thick] (0,0)--(10,0) node[above left]{$L$};
 \node[dot,label=below:{$p$}] at (.58,0){};
 \node[dot] (q) at (1.59,0){};
 \node at (1.3,-.3){$q$};
 \node[gate,label=above:{$F$}] at (2.6,0){};
 \node[gate] at (2.2,1.19){};
 \node[gate] at (4.3,-2.5){};
 \draw[rounded corners=8pt,dashed] (1.2,-3.06) rectangle (6.89,.88);
 \node at (6.5,.58){\small$Z^*$};
 \draw[dotted] (1.2,-1.47)--(6.89,-1.47);
 \node at (6.5,-1.18){\small$A^*$};
 \node[odot] (o1) at (2.49,-.7){};
 \node[dot] (d1) at (3.59,-.51){};
 \node[dot] (d2) at (3.49,-1.15){};
 \node[odot] (o2) at (4.7,-.7){};
 \node[dot] (d3) at (2.39,-2.11){};
 \node[dot] (d4) at (3.39,-2.3){};
 \draw[->] (q)--(o1);
 \draw[->] (o1)--(d1);
 \draw[->] (o1)--(d2);
 \draw[->] (d1)--(o2);
 \draw[->] (d2)--(o2);
 \draw[->] (o2) to[bend right=10] (7.47,-.08);
 \draw[->] (o1)--(d3);
 \draw[->] (d2)--(d4);
 \draw[->] (d3)--(d4);
 \draw[->,thick,gray] (d4)--(5.6,-2.3) node[above left,black]{\scriptsize$L'$};
 \node at (6.08,-2.54){\small$B=Z'$};
 \draw[decorate,decoration={brace,mirror,amplitude=4pt}] (7.41,-.31)--(9.78,-.31);
 \node at (8.6,-.75){\small tail $T$};
\end{tikzpicture}

\smallskip
{\small The part $A^*$ of $Z^*$ can still reach the tail $T$ without
passing the gate. All paths from $q$ to $T$ through $A$ have the
same type, so $A^*$ is two-coloured by the parity of the distance
from $q$. The part $B$ cannot reach $T$. The paradox, and a new
ladder ray $L'$, are in $B$.}
\end{center}

\begin{lemma}[Gate]
\label{n6.1}
Let $Z$ be a finite-port region of a noc graph $G$, let $L$ be a ray, let $q\in
Z$ be an essential vertex of $G$ that does not parity-dominate $L$,
and let $A\not\ni q$ be a finite set. Then there is a finite-port
region $Z'\subseteq Z\setminus A$ of $G$ such that
\begin{enumerate}
\item every vertex of $Z'$ is reachable from $q$ by a path in
$G[Z]-A$;
\item no ray lying in $Z'$ feeds $L$;
\item some ladder ray of $G$ lies in $Z'$.
\end{enumerate}
\end{lemma}
\begin{proof}
By the parity cut criterion (Lemma \ref{n3.3}) there are a finite
set $F\not\ni q$ with $A\subseteq F$ and a tail $T$ of $L$ with
$q\notin T$ such that all paths of positive length in $G-F$ from $q$
to $T$ have the same type with respect to $L$. By the parity
principle the same is true of all walks in $G-F$ from $q$ to $T$:
such a walk contains a path with the same ends and the same parity.

Let $Z^*$ be the cone of vertices reachable from $q$ in $G[Z]-F$. By
Lemma \ref{n4.12}(2), $G[Z^*]$ is not kernel-perfect. Let
\[
 A^*=\{x\in Z^*:\ \text{some path in $G-F$ leads from $x$ to $T$}\},
 \qquad B=Z^*\setminus A^* .
\]

\emph{$A^*$ is bipartite.} Let $x\in A^*$, and let $P,P'$ be walks
from $q$ to $x$ in $G[Z]-F$. Prolong both by the same path in $G-F$
from $x$ to a vertex of $T$. The two walks obtained have the same
type and the same last vertex, so the same parity, and so do $P$ and
$P'$. Colour $x$ by this common parity. If $x\to y$ is an edge with
$x,y\in A^*$, a walk from $q$ to $x$ followed by this edge is a walk
from $q$ to $y$, so $x$ and $y$ have different colours. Hence
$G[A^*]$ is bipartite, and kernel-perfect by Lemma \ref{n2.6}.

\emph{$B$ is closed.} Let $x\to y$ be an edge with $x\in B$ and
$y\in Z^*$. If $y$ had a path to $T$ in $G-F$, so would $x$. So $B$
is successor-closed in $G[Z^*]$.

\emph{The paradox is in $B$.} If $G[B]$ were kernel-perfect, then
$G[Z^*]$ would be kernel-perfect by Lemma \ref{n2.5}. So $G[B]$ is
not kernel-perfect. Put $Z'=B$. An edge that leaves $B$ leaves
$Z^*$, so $\bd B\subseteq\bd Z^*$ is finite. Property (1) holds for
all of $Z^*$, and (3) holds by Proposition \ref{n5.5}. For (2),
suppose a ray $S$ lying in $Z'$ feeds $L$. By Lemma \ref{n3.7}(1)
there is a path in $G-F$ from a vertex of $S$ to $T$. That vertex
belongs to $A^*$, not to $B$.
\end{proof}

Each hypothesis is used once: the absence of a mixed fan produces
the gate and makes $A^*$ bipartite; essentiality of $q$ keeps the
paradox inside $Z^*$, hence inside $B$; and the cover lemma,
through Proposition \ref{n5.5}, turns that paradox into a new ladder
ray. If $q$ does not dominate $L$ at all, $A^*$ can be made empty;
this is the gate of a safe digraph. The set $A$ is there to keep later
constructions away from vertices already used.

\subsection{The diagonal comb and the main result}\label{s6.2}
By the following lemma, a parity-safe $\Gamma$ having some ladder end has uncountably many.
\begin{lemma}[Diagonal comb]
\label{n6.3}
If $\Gamma$ is parity-safe, $Z$ a finite-port region of $\Gamma$ in which some ladder ray lies, and
$e_0,e_1,e_2,\dots$ a sequence of (not necessarily distinct) ladder ends of $\Gamma$, then
there are finite-port regions $Z=Z_0\supseteq Z_1\supseteq Z_2\supseteq\cdots$ of $\Gamma$ and
a ladder ray $R$ lying in $Z$ such that for every $i$

1. $R$ has a tail that lies in $Z_i$, and 

2. no ladder ray lying in $Z_{i+1}$ belongs to $e_i$.
\end{lemma}

The picture is the one in Section \ref{s1.3}, given with more details below.
 The ray $R$ is the
\emph{spine} of the comb and the rays $L_i$ are its \emph{teeth}.

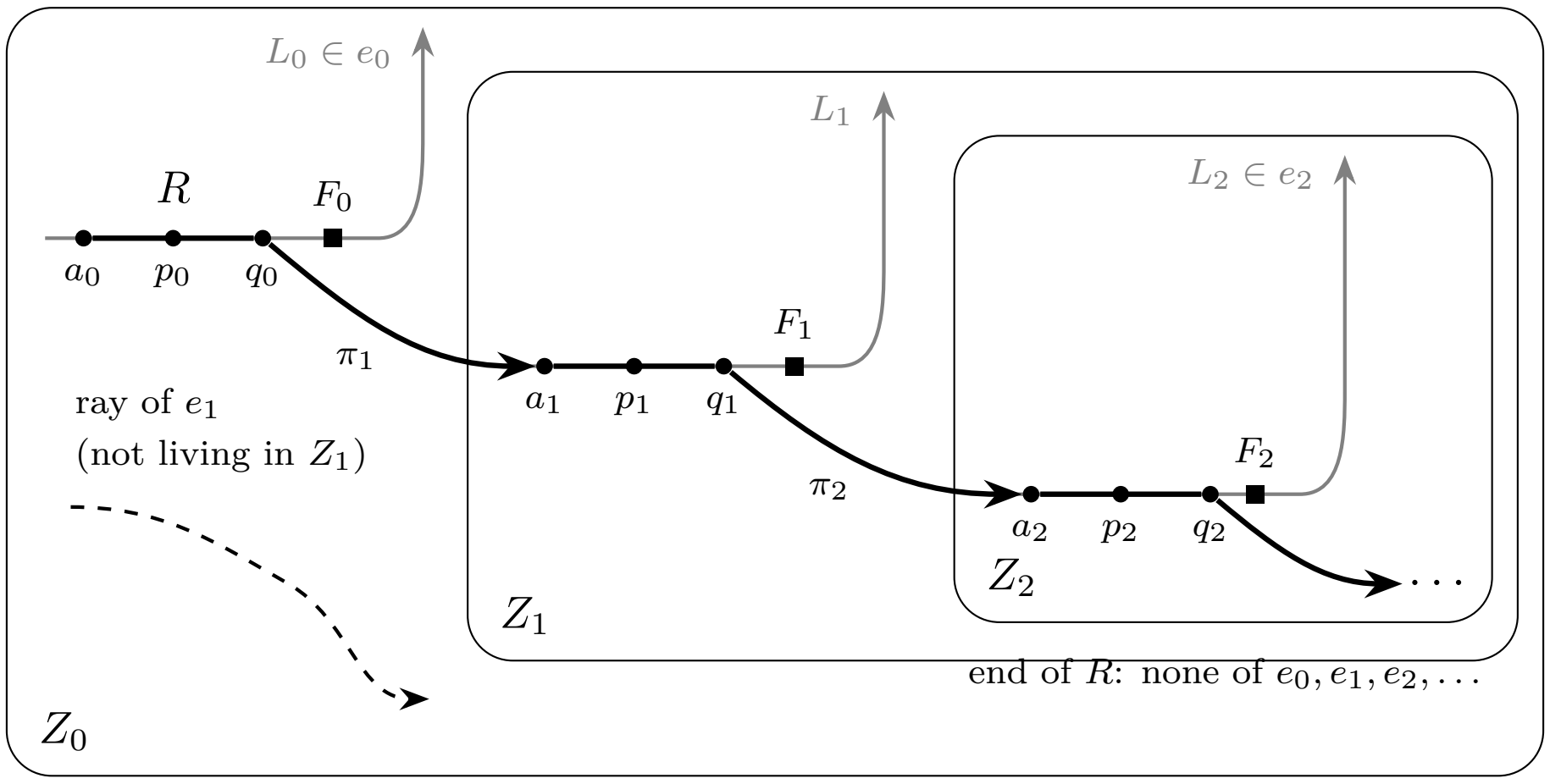
\begin{figure}[htbp]
\centering
\begin{tikzpicture}
 % nested regions
 \draw[rounded corners=10pt] (0,-.6) rectangle (12,5.4);
 \draw[rounded corners=10pt] (3.6,.3) rectangle (11.8,4.9);
 \draw[rounded corners=10pt] (7.4,.6) rectangle (11.6,4.4);
 \node at (.45,-.25){\small$Z_0$};
 \node at (4.05,.65){\small$Z_1$};
 \node at (7.85,.95){\small$Z_2$};

 % tooth L_0, a ray of e_0
 \draw[->,gray,thick] (.3,3.6)--(2.9,3.6) to[out=0,in=-90] (3.25,4.4)--(3.25,5.25);
 \node[gate,label=above:{\scriptsize$F_0$}] at (2.55,3.6){};
 \node[dot,label=below:{\scriptsize$a_0$}] (a0) at (.6,3.6){};
 \node[dot,label=below:{\scriptsize$p_0$}] at (1.3,3.6){};
 \node[dot,label=below:{\scriptsize$q_0$}] (q0) at (2,3.6){};
 \node[gray,anchor=east] at (3.15,5.05){\scriptsize$L_0\in e_0$};

 % e_1 does not live in Z_1: one of its ladder rays, lying outside Z_1
 \draw[->,thick,dashed] (.5,1.5) to[out=0,in=150] (2.2,.9) to[out=-30,in=180] (3.3,0);
 \node[anchor=west] at (.4,2.3){\scriptsize ray of $e_1$};
 \node[anchor=west] at (.4,1.9){\scriptsize (not living in $Z_1$)};

 % tooth L_1, an arbitrary ladder ray in Z_1
 \draw[->,gray,thick] (3.9,2.6)--(6.5,2.6) to[out=0,in=-90] (6.85,3.4)--(6.85,4.75);
 \node[gate,label=above:{\scriptsize$F_1$}] at (6.15,2.6){};
 \node[dot,label=below:{\scriptsize$a_1$}] (a1) at (4.2,2.6){};
 \node[dot,label=below:{\scriptsize$p_1$}] at (4.9,2.6){};
 \node[dot,label=below:{\scriptsize$q_1$}] (q1) at (5.6,2.6){};
 \node[gray,anchor=east] at (6.75,4.6){\scriptsize$L_1$};

 % tooth L_2, a ray of e_2
 \draw[->,gray,thick] (7.7,1.6)--(10.1,1.6) to[out=0,in=-90] (10.45,2.4)--(10.45,4.25);
 \node[gate,label=above:{\scriptsize$F_2$}] at (9.75,1.6){};
 \node[dot,label=below:{\scriptsize$a_2$}] (a2) at (8,1.6){};
 \node[dot,label=below:{\scriptsize$p_2$}] at (8.7,1.6){};
 \node[dot,label=below:{\scriptsize$q_2$}] (q2) at (9.4,1.6){};
 \node[gray,anchor=east] at (10.35,4.1){\scriptsize$L_2\in e_2$};

 % spine
 \draw[very thick] (a0)--(q0);
 \draw[->,very thick] (q0) to[out=-40,in=180] node[below left=-1pt]{\scriptsize$\pi_1$} (a1);
 \draw[very thick] (a1)--(q1);
 \draw[->,very thick] (q1) to[out=-40,in=180] node[below left=-1pt]{\scriptsize$\pi_2$} (a2);
 \draw[very thick] (a2)--(q2);
 \draw[->,very thick] (q2) to[out=-40,in=180] (10.9,.9);
 \node at (11.2,.9){$\cdots$};
 \node at (1.3,4.0){\small$R$};
 \node[anchor=north east] at (11.7,.45){\scriptsize end of $R$: none of $e_0,e_1,e_2,\dots$};
\end{tikzpicture}

\smallskip {\small 
  At stage $i$, if $e_i$ lives in $Z_i$ the tooth $L_i$ is a ray of $e_i$ (here $e_0$ and
  $e_2$); otherwise $L_i$ is any ladder ray lying in $Z_i$ (here $e_1$, whose rays lie
  outside $Z_1$). The spine $R$ (thick) follows $L_i$ only from $a_i$ to $q_i$; behind the
  gate $F_i$ the tail of $L_i$ turns away and stays outside $Z_{i+1}$. No ladder ray lying
  in $Z_{i+1}$ belongs to $e_i$: in the first case it would feed $L_i$, which the gate
  excludes; in the second $e_i$ does not live even in $Z_i\supseteq Z_{i+1}$. A tail of
  $R$ lies in every $Z_{i+1}$, so the end of $R$ is none of the $e_i$.}
\end{figure}

\begin{proof}
  Being parity-safe, $\Gamma$ has no odd cycles, so the results stated above for noc
  graphs apply to it.

  (1). We build $R$ as an increasing union of paths $R_0\subseteq R_1\subseteq\cdots$,
  where $R_i$ ends at a vertex $q_i$. Let $Z_0=Z$. Suppose $Z_i$ is given and some ladder
  ray lies in it, and, if $i>0$, that $R_{i-1}$ is given, that $Z_i$ contains no vertex of
  $R_{i-1}$ other than possibly $q_{i-1}$, and that every vertex of $Z_i$ is reachable
  from $q_{i-1}$ by a path in $\Gamma[Z_{i-1}]$ that meets $R_{i-1}$ only in $q_{i-1}$.
\begin{itemize}
\item \emph{The tooth.} If $e_i$ lives in $Z_i$, let $L_i$ be a ladder ray of $e_i$ lying
  in $Z_i$. Otherwise let $L_i$ be any ladder ray lying in $Z_i$.
\item \emph{The entry.} If $i=0$, let $a_0$ be the first vertex of $L_0$ ($p_0$ and $q_0$
  are determined as in the next step). If $i>0$ and $q_{i-1}$ lies on $L_i$, let
  $a_i=q_{i-1}$.  Otherwise take a path $\pi_i$ as above from $q_{i-1}$ to a vertex of
  $L_i$, stop it at its first vertex $a_i$ on $L_i$, and append it to $R_{i-1}$.
\item \emph{Along the tooth.} $L_i$ has infinitely many checkpoints beyond $a_i$, all
  essential, and by parity-safety only finitely many of them parity-dominate $L_i$. Choose
  consecutive checkpoints $p_i\prec q_i$ beyond $a_i$ such that $q_i$ does not
  parity-dominate $L_i$. Append the segment of $L_i$ from $a_i$ to $q_i$. The result is
  $R_i$. It is a path: the segment lies in $Z_i$, which avoids $R_{i-1}$ except for
  $q_{i-1}$, and it meets $\pi_i$ only in $a_i$.
\item \emph{The gate.} Apply Lemma \ref{n6.1} to $Z_i$, $L_i$, $q_i$ and
  $A=V(R_i)\setminus\{q_i\}$. This gives $Z_{i+1}\subseteq Z_i$, which satisfies the
  assumptions made above for the next step.
\item \emph{The spine.} $R=\bigcup_i R_i$ is a ray, since $R_i$ properly extends
  $R_{i-1}$. It passes through $p_0,q_0,p_1,q_1,\dots$ in this order, the $p_i$ are
  essential and $p_i\prec q_i$; so $R$ is a ladder ray by Lemma \ref{n5.3}(2). It lies in
  $Z_0$. From $a_{i+1}$ on, $R$ consists of segments of $L_j$ with $j>i$, which lie in
  $Z_j\subseteq Z_{i+1}$, and of paths $\pi_j$ with $j>i+1$, which lie in
  $Z_{j-1}\subseteq Z_{i+1}$. This is (1).
\end{itemize}
(2). Suppose a ladder ray $S$ lying in $Z_{i+1}$ belongs to $e_i$. As
$Z_{i+1}\subseteq Z_i$, the end $e_i$ lives in $Z_i$, so $L_i$ was chosen in $e_i$, and
$S$ feeds $L_i$. This contradicts Lemma \ref{n6.1}(2).
\end{proof}
This gives the following theorem.

\begin{theorem}\label{n6.4}
Let $\Gamma$ be a parity-safe digraph that is not kernel-perfect, and $Z$
be a finite-port region of $\Gamma$ that contains an essential vertex.
Then uncountably many ladder ends live in $Z$. 
\end{theorem}
\begin{proof}
Being parity-safe, $\Gamma$ has no odd cycles. $\Gamma[Z]$ is not kernel-perfect, so a ladder ray lies in $Z$ by
Proposition \ref{n5.5}, and its end lives in $Z$. Suppose the ladder
ends living in $Z$ could be listed as $e_0,e_1,e_2,\dots$, with
repetitions if there are finitely many. Take the regions $Z_i$ and
the ladder ray $R$ of Lemma \ref{n6.3}. Since $R$ lies in $Z$, its
end is a ladder end living in $Z$, so it is $e_i$ for some $i$. A
tail of $R$ lies in $Z_{i+1}$. It is a ladder ray and belongs to the
end of $R$, by Lemma \ref{n5.3}(1). So a ladder ray lying in
$Z_{i+1}$ belongs to $e_i$, contrary to Lemma \ref{n6.3}(2).

For the last sentence: $\Gamma$ is a noc graph that is not kernel-perfect, so it has an
essential vertex, by Lemma \ref{n4.12}(3) and Corollary \ref{n2.4},
and $V(\Gamma)$ is a finite-port region containing it.
\end{proof}
The main result follows now easily. 

\begin{restated}{n1.2}
A parity-safe digraph with at most countably many ends is
kernel-perfect. In particular, a safe digraph with at most countably
many ends is kernel-perfect.
\end{restated}
\begin{proof}
If a parity-safe $\Gamma$ has countably many ends, then it has no essential vertex by Theorem \ref{n6.4}, hence no core, by Lemma \ref{n4.12}, and hence is kernel-perfect by Corollary \ref{n2.4}.
Safe digraphs are parity-safe.
\end{proof}

\section{Lean definitions and Theorem \ref{n1.2}}\label{s8}
{\small{\tt 
variable $\{$V : Type*$\}$ (E : V → V → Prop)
\\
def IsCycle (v : V) (p : List V) : Prop :=

p $\not=$ [] $\land$ List.IsChain E (v :: p) $\land$ 

p.getLast? = some v $\land$ p.Nodup
\\
def NoOddCycle : Prop :=

  $\forall$ (v : V) (p : List V), IsCycle E v p → Even p.length
\\
def IsKernelOf (X K : Set V) : Prop :=

  K $\subseteq$ X $\land$

 ($\forall$ u $\in$ K, $\forall$ w $\in$ K, ¬ E u w) $\land$

  ($\forall$ v $\in$ X, v $\not\in$ K → $\exists$ w $\in$ K, E v w)
\\
def KernelPerfect : Prop :=

  $\forall$ X : Set V, $\exists$ K : Set V, IsKernelOf E X K
\\
def IsRay (r : Nat → V) : Prop :=

  Function.Injective r $\land$ 

$\forall$ n : Nat, E (r n) (r (n + 1))
\\
def IsFan (v : V) (R : Nat → V) (P : Nat → List V) (k : Nat → Nat) : Prop :=

  ($\forall$ n : Nat,  P n $\not=$ [] $\land$ 

List.IsChain E (v :: P n) $\land$ 

(v :: P n).Nodup $\land$

   (P n).getLast? = some (R (k n))) $\land$

  ($\forall$ m n : Nat, m $\not=$ n → List.Disjoint (P m) (P n))
% \\
% def Dominates (v : V) (R : Nat → V) : Prop :=
%
%   $\exists$ (P : Nat → List V) (k : Nat → Nat), IsFan E v R P k
\\
def ParityDominates (v : V) (R : Nat → V) : Prop :=
  
  $\exists$ (P : Nat → List V) (k : Nat → Nat), IsFan E v R P k $\land$

  $\{$n : Nat | Even ((P n).length + k n)$\}$.Infinite $\land$

    $\{$n : Nat | Odd ((P n).length + k n)$\}$.Infinite
% \\
% def Safe : Prop :=
%
%   NoOddCycle E $\land$
%   $\forall$ R : Nat → V, IsRay E R → ¬ $\{$n : Nat | Dominates E (R n) R$\}$.Infinite
\\
def IsParitySafe : Prop :=

  NoOddCycle E $\land$
  $\forall$ R : Nat → V, IsRay E R → ¬ $\{$n : Nat | ParityDominates E (R n) R$\}$.Infinite
\\
def ReachAvoiding (F : Finset V) (u w : V) : Prop :=

   u $\not\in$ F $\land$ ReflTransGen (fun a b => E a b $\land$ b $\not\in$ F) u w
\\
def Feeds (R R' : Nat → V) : Prop :=

   $\forall$ F : Finset V, $\exists$ m n : Nat, ReachAvoiding E F (R m) (R' n)
\\
def endOf (R : Nat → V) : Set (Nat → V) :=
  $\{$R' | IsRay E R' $\land$ Feeds E R R' $\land$ Feeds E R' R$\}$
\\
def Ends : Set (Set (Nat → V)) :=
  $\{$S | $\exists$ R : Nat → V, IsRay E R $\land$ S = endOf E R$\}$
\\
\\
/-- **Theorem 1.2.** A parity-safe digraph with countably many ends is kernel-perfect. -/
\\
theorem kernelPerfect\_of\_paritySafe\_of\_countable\_ends

  (hS : IsParitySafe E) (hEnds : (Ends E).Countable) :  KernelPerfect E := by ...
}}

\end{document}